\documentclass{article}

\usepackage{PRIMEarxiv}

\usepackage[utf8]{inputenc} % allow utf-8 input
\usepackage[T1]{fontenc}    % use 8-bit T1 fonts
\usepackage{hyperref}       % hyperlinks
\usepackage{url}            % simple URL typesetting
\usepackage{booktabs}       % professional-quality tables
\usepackage{amsfonts}       % blackboard math symbols
\usepackage{nicefrac}       % compact symbols for 1/2, etc.
\usepackage{microtype}      % microtypography
\usepackage{lipsum}
\usepackage{enumitem} % for custom enumerate labels
\usepackage{fancyhdr}       % header
\usepackage{graphicx}       % graphics
\graphicspath{{media/}}     % organize your images and other figures under media/ folder
\usepackage{siunitx}
\usepackage{amsfonts}
\usepackage{amsmath}
\usepackage{amsthm}
\usepackage{enumitem} 
\usepackage{url}

\newtheorem{theorem}{Theorem}

\newtheorem{definition}{Definition}

\newtheorem{lemma}{Lemma}

\newtheorem{remark}{Remark}

\title{ A Complete Characterization Theorem for Fuzzy Differentiability on Time Scales
}

\author{
  Funda Raziye MERT \\
  Department of Software Engineering \\
  Adana Alparslan Türkeş Science and Technology University \\
  Adana, TÜRKİYE\\
  \texttt{rmert@atu.edu.tr} \\
   \And
  Selami BAYEĞ \\
  Department of Industrial Engineering\\
  University of Turkish Aeronautical Association \\
  Ankara, TÜRKİYE\\
  \texttt{sbayeg@thk.edu.tr} \\
   \AND
}

\begin{document}
\maketitle

\begin{abstract}

This paper investigates the generalized Hukuhara differentiability of fuzzy number-valued functions on arbitrary time scales using delta calculus. By carefully examining and improving existing results, we develop a unified and complete characterization theorem that covers a wide range of differentiability behaviors, including some cases that were previously missed. Our approach addresses important limitations and redundancies in earlier work, providing a clearer and more flexible understanding of fuzzy differentiability.

\end{abstract}

	\section{Introduction}

\qquad In mathematical modeling and analysis, time plays a fundamental role in describing the evolution of dynamic systems. Traditional approaches often treat time as either continuous, using differential equations, or discrete, using difference equations. However, many real-world systems exhibit behavior that does not fit neatly into either category, involving both continuous changes and discrete events. The theory of time scales, first introduced by Stefan Hilger in the late 1980s, offers a unified framework that accommodates both discrete and continuous time within a single mathematical structure ~\cite{Hilger1}.

Time scale calculus enables the study of dynamic equations on arbitrary closed subsets of the real numbers, known as time scales, thus generalizing classical results and allowing for broader applicability. A central concept within this theory is delta calculus, which extends the derivative to time scales, offering a coherent way to describe the rate of change regardless of whether the time domain is continuous, discrete, or a combination of both. This approach has proven particularly useful in fields such as control theory, economics, biology, and engineering, where systems often operate across mixed temporal regimes. For detailed information on time scales, see the monographs~\cite{Peterson, Peterson1}.

Fuzzy set theory, introduced by Lotfi Zadeh in 1965, extends the framework of classical set theory by allowing partial membership, where elements can belong to a set to varying degrees between 0 and 1~\cite{Zadeh}. Unlike traditional sets, which impose a strict binary classification, fuzzy sets accommodate uncertainty and gradual transitions. This flexibility has made fuzzy theory a powerful approach in fields that involve imprecise or vague information, such as intelligent systems, control engineering, decision analysis, and pattern recognition~\cite{Klir, Gomes, Bedebook}. By reflecting the ambiguity inherent in many real-world situations, fuzzy set theory supports more realistic and adaptive models. 

The generalized Hukuhara difference represents a refined extension of the classical Hukuhara difference, designed to address limitations encountered in fuzzy set theory when comparing fuzzy numbers. This generalization provides a more adaptable and comprehensive mechanism for handling operations involving fuzzy quantities, especially within the study of dynamic systems and fuzzy differential equations. Its utility becomes particularly evident when working with models that span both continuous and discrete time settings. The development of the generalized Hukuhara delta derivative further enhances the standard delta derivative by extending its applicability to fuzzy-valued functions on arbitrary time scales. When combined with time scale calculus, this approach offers a unified and powerful methodology for analyzing fuzzy dynamic systems, particularly in scenarios where traditional difference or differential methods fall short, such as in hybrid or uncertain environments.

In this paper, we investigate the generalized Hukuhara differentiability of fuzzy number-valued functions on arbitrary time scales using the tools of delta calculus. We build on and refine the foundational results from \cite{Cano, qiu, longo} to present a complete and unified characterization theorem. Chalco-Cano et al. \cite{Cano} initially characterized gH-differentiability based on the differentiability of endpoint functions. However, Qiu \cite{qiu} demonstrated through a counterexample that this approach is not always valid for interval-valued functions, and introduced a more general framework that includes atypical cases—such as functions that are gH-differentiable at a point but discontinuous in any punctured neighborhood. Longo et al. \cite{longo} later extended these ideas to fuzzy number-valued functions. Our work further develops this theory by addressing the limitations and redundancies in Longo’s formulation, resulting in a more robust and broadly applicable framework for fuzzy differential calculus across discrete, continuous, and hybrid time domains.

This paper is organized into three sections. The preliminaries section lays the groundwork by presenting essential definitions, foundational lemmas, and critical theorems that support the subsequent theoretical developments. In the main results section, we introduce a comprehensive characterization theorem for the delta generalized Hukuhara differentiability of fuzzy number-valued functions on arbitrary time scales. Finally, the conclusion section offers a concise summary of our findings and highlights their significance.

%In this paper, we have modified the definition of the generalized Hukuhara delta  derivative introduced by Fard et al. \cite{Fard} and extensively examined its properties, with the goal of addressing gaps in the existing literature while substantially broadening and enhancing the scope of previous results. We provide several characterizations of generalized Hukuhara delta  differentiable fuzzy functions on time scales, based on the delta differentiability of their endpoint functions. These characterizations serve as valuable tools for the calculus of derivatives of fuzzy functions on time scales. Our findings extend the results previously obtained in \cite{Cano} for fuzzy number valued functions in the real case to time scale case. We have also extended the product rule, which was proven in \cite{TZ} for interval valued functions in the real case to fuzzy number valued functions on time scales.

%The paper provides useful tools for the calculus of derivatives of fuzzy functions because it gives expressions of the
%level sets of the derivatives in terms of the derivatives of the endpoint functions. These expressions have been given
%for several types of differentiability, specifically, for generalized Hukuhara differentiable fuzzy functions. 

\section{Preliminaries}
This section introduces fundamental definitions, preliminary lemmas, and key theorems that will underpin the theoretical development in the following sections.

The definitions presented here are primarily based on the framework established by Bohner et al. \cite{Peterson}.

A time scale \(\mathbb{T}\) is any closed subset of the real numbers \(\mathbb{R}\). The forward and backward jump operators, \(\sigma\) and \(\rho\), map \(\mathbb{T}\) into itself and are defined by
\[
\sigma(t) = \inf\{ s \in \mathbb{T} : s > t \}, \quad \rho(t) = \sup\{ s \in \mathbb{T} : s < t \},
\]
with the conventions \(\inf \emptyset := \sup \mathbb{T}\) and \(\sup \emptyset := \inf \mathbb{T}\). A point \( t \in \mathbb{T} \) is classified as left-dense if \( t > \inf \mathbb{T} \) and \(\rho(t) = t\); left-scattered if \(\rho(t) < t\); right-dense if \( t < \sup \mathbb{T} \) and \(\sigma(t) = t\); and right-scattered if \(\sigma(t) > t\). If \(\mathbb{T}\) has a maximum point \( M \) that is left-scattered, then define \(\mathbb{T}^{\kappa} := \mathbb{T} \setminus \{M\}\); otherwise, \(\mathbb{T}^{\kappa} := \mathbb{T}\). The forward graininess function \(\mu : \mathbb{T} \to [0, \infty)\) is given by
\[
\mu(t) = \sigma(t) - t,
\]
representing the distance to the next point in the time scale.

For a function $f : \mathbb{T} \to \mathbb{R}$ and a point $t \in \mathbb{T}^{\kappa}$, the delta derivative of $f$ at $t$, denoted by $\Delta f(t)$, is defined as the number (if it exists) such that for any $\varepsilon > 0$, there is a neighborhood $U$ of $t$ satisfying
$$|[f(\sigma(t)) - f(s)] - \Delta f(t)[\sigma(t) - s]| \leq \varepsilon |\sigma(t) - s|$$
for all $s \in U$.
\begin{theorem}\label{nabladerrule}\cite{Peterson}
Let \( f: \mathbb{T} \to \mathbb{R} \) be a function, and let \( t \in \mathbb{T}^\kappa \) be a point. The following statements hold:
\begin{enumerate}[label=(\roman*)]
\item If \( f \) is $\Delta$-differentiable at \( t \), then \( f \) is continuous at \( t \).
\item If \( f \) is continuous at \( t \) and \( t \) is right-scattered, then \( f \) is $\Delta$-differentiable at \( t \) with
\begin{equation*}
\Delta f(t) = \frac{f(\sigma(t)) - f(t)}{\mu(t)}.
\end{equation*}
\item If \( t \) is right-dense, then \( f \) is $\Delta$-differentiable at \( t \) if and only if the limit
$$\lim_{s \to t} \frac{f(s) - f(t)}{s - t}$$
exists. In this case, we have
$$\Delta f(t) = \lim_{s \to t} \frac{f(s) - f(t)}{s - t}.$$
\item If \( f \) is $\Delta$-differentiable at \( t \), then
\begin{equation}
f(\sigma(t)) = f(t) + \mu(t) \Delta f(t).
\end{equation}
\end{enumerate}
\end{theorem}

\begin{definition}\cite{Zadeh}
A fuzzy set \( u \) on a universe \( U \) is characterized by a membership function \( u : U \to [0,1] \), where \( u(x) \) indicates the membership grade of element \( x \) in \( u \).
\end{definition}

\begin{definition}\cite{Negoita}
Given a fuzzy set \( u: U \to [0,1] \), its \(\alpha\)-level set for \( \alpha \in (0,1] \) is defined as
\[
[u]_\alpha = \{ x \in U : u(x) \geq \alpha \},
\]
while the support, or \( 0 \)-level set, is given by
\[
[u]_0 = \overline{\{ x \in U : u(x) > 0 \}},
\]
where the bar denotes the closure.
\end{definition}

\begin{definition}\cite{Bedebook}
A fuzzy number \( u: \mathbb{R} \to [0,1] \) satisfies the following:
\begin{enumerate}
    \item \textbf{Normality:} There exists \( x_0 \in \mathbb{R} \) with \( u(x_0) = 1 \).
    \item \textbf{Fuzzy Convexity:} For any \( t \in [0,1] \) and \( x,y \in \mathbb{R} \),
    \[
    u(tx + (1-t)y) \geq \min\{u(x), u(y)\}.
    \]
    \item \textbf{Upper semicontinuity:} For each \( x_0 \in \mathbb{R} \) and any \( \varepsilon > 0 \), there exists \( \delta > 0 \) such that
\[
u(x)-u(x_0)< \varepsilon \quad \text{whenever} \quad |x - x_0| < \delta.
\]
    \item \textbf{Compact Support:} The support of \( u \) is compact.
\end{enumerate}
\end{definition}

The set of all fuzzy numbers is denoted by \( \mathbb{R}_{\mathcal{F}} \).

\begin{theorem}\cite{Bedebook}\label{Bede1}
Let \( u^-, u^+ : [0,1] \to \mathbb{R} \) be bounded functions where \( u^- \) is non-decreasing, left-continuous on \((0,1]\), and right-continuous at 0; similarly, \( u^+ \) is non-increasing, left-continuous on \((0,1]\), and right-continuous at 0. Assume \( u^-_1 \leq u^+_1 \). Then there is a unique fuzzy number \( u \in \mathbb{R}_{\mathcal{F}} \) with \(\alpha\)-level sets \( [u]_\alpha = [u^-_\alpha, u^+_\alpha] \).
\end{theorem}

\begin{definition}\cite{Stef}
The generalized Hukuhara difference (\(gH\)-difference) between fuzzy numbers \( u \) and \( v \) is the fuzzy number \( \omega \), if it exists, such that
\[
u \ominus_{gH} v = \omega \quad \text{iff} \quad u = v + \omega \quad \text{or} \quad v = u + (-1)\omega.
\]
\end{definition}

\begin{remark}\cite{Stef,Bede}
For the existence of \( \omega = u \ominus_{gH} v \), one of the following must hold:

\textbf{Case (i):}  
\begin{itemize}
    \item \( \omega^-_\alpha = u^-_\alpha - v^-_\alpha \) and \( \omega^+_\alpha = u^+_\alpha - v^+_\alpha \) for all \(\alpha \in [0,1]\),
    \item \( \omega^-_\alpha \) is non-decreasing and \( \omega^+_\alpha \) is non-increasing with \( \omega^-_\alpha \leq \omega^+_\alpha \) for all \(\alpha \in [0,1]\).
\end{itemize}

\textbf{Case (ii):}  
\begin{itemize}
    \item \( \omega^-_\alpha = u^+_\alpha - v^+_\alpha \) and \( \omega^+_\alpha = u^-_\alpha - v^-_\alpha \) for all \(\alpha \in [0,1]\),
    \item \( \omega^-_\alpha \) is non-decreasing and \( \omega^+_\alpha \) is non-increasing with \( \omega^-_\alpha \leq \omega^+_\alpha \) for all \(\alpha \in [0,1]\).
\end{itemize}
\end{remark}

\begin{theorem}\cite{Stef, Bede}\label{prop2}
If \( \omega = u \ominus_{gH} v \) exists for \( u,v \in \mathbb{R}_{\mathcal{F} } \), then the \(\alpha\)-level sets of \(\omega\) satisfy
\[
[\omega]_\alpha = [\omega^-_\alpha, \omega^+_\alpha],
\]
where
\[
\omega^-_\alpha = \min\{ u^-_\alpha - v^-_\alpha, u^+_\alpha - v^+_\alpha \}, \quad \omega^+_\alpha = \max\{ u^-_\alpha - v^-_\alpha, u^+_\alpha - v^+_\alpha \}.
\]
\end{theorem}

\begin{lemma}\cite{Cano}
Given fuzzy numbers \( u,v \) with an existing \( u \ominus_{gH} v \), then for any scalar \( a \in \mathbb{R} \), the difference \( a u \ominus_{gH} a v \) also exists and satisfies
\[
a (u \ominus_{gH} v) = a u \ominus_{gH} a v.
\]
\end{lemma}

The Hausdorff metric between fuzzy numbers, essential in applications such as clustering, pattern recognition, and decision-making, is defined as follows. This metric generalizes the classical Hausdorff distance to handle the uncertainty intrinsic to fuzzy numbers, thereby providing a reliable measure of similarity between fuzzy quantities.

\begin{definition}\cite{Diamond2} \label{def12}
	The Hausdorff metric $D:\mathbb{R}_{\mathcal{F} }\times \mathbb{R}_{\mathcal{F} }
    \rightarrow 
	[0,\infty)$ is defined by
    \begin{equation*}
    \begin{split}
	D(u,v)&=\sup_{\alpha\in \lbrack 0,1]}D_H([u]_{\alpha},[v]_{\alpha})\\
    &=\sup_{\alpha\in \lbrack 0,1]}\max \left\{ \left\vert
	u_{\alpha}^{-}-v_{\alpha}^{-}\right\vert ,\left\vert u_{\alpha}^{+}-v_{\alpha}^{+}\right\vert
	\right\}, 
    \end{split}
    \end{equation*}
	where $[u]_{\alpha}=[u_{\alpha}^{-},u_{\alpha}^{+}]$, $[v]_{\alpha}=[v_{\alpha}^{-},v_{\alpha}^{+}]$.
\end{definition}

%The Hausdorff metric provides a way to measure the distance between two fuzzy sets by considering their level sets. This metric allows researchers to compare the similarity or dissimilarity of fuzzy sets in a rigorous mathematical way. Specifically, it can be used to quantify how far apart two fuzzy sets are based on their support and their membership functions.

\begin{theorem}\cite{Diamond2} \label{prop3} Let $u, v, w, e\in \mathbb{R}_{\cal{F}}$ and $k\in\mathbb{R}$. The Hausdorff metric satisfies the followings:
	\begin{enumerate}
		\item [(i)] $D\left( u+w,v+w\right) =D\left(u,v\right)$

  	\item [(ii)] $D(ku,kv) =\left\vert k\right\vert
		D\left(u,v\right)$

		\item [(iii)] $D\left(u+v,w+e\right) \leq D\left(
		u,w\right) +D\left(v,e\right). $
	\end{enumerate}
    %and $(\mathbb{R}_{\cal{F}},D_{\infty })$ is a complete metric space.
	
\end{theorem}

\begin{definition}\cite{Zakon}\label{def1}
Let $\{x_n\}_{n=1}^{\infty}$ be a sequence of real numbers, and let $c \in [-\infty,\infty]$. The point  $c$ is called a \textit{cluster point}
 (or accumulation point) of the sequence if every neighborhood of $c$ contains infinitely many terms of the sequence.
\end{definition}

It is a well-known result that  $c \in [-\infty,\infty]$ is a cluster point of the sequence \( \{x_n\}_{n=1}^{\infty} \) if and only if there exists a subsequence \( \{x_{n_k}\}_{k=1}^{\infty} \) such that
$
\lim_{k \to \infty} x_{n_k} = c.$
\begin{definition}
Let \( f : \mathbb{T} \to \mathbb{R} \) be a real valued function, and let \( t_0 \in \mathbb{T} \) be a point.  
A number \( c \in [-\infty, \infty] \) is called a \emph{right-sided cluster point} of \( f \) at \( t_0 \) if there exists a sequence of positive real numbers \( \{h_n\}_{n=1}^{\infty} \) such that:
\begin{enumerate}[label=(\roman*)]
    \item \( \displaystyle\lim_{n \to \infty} h_n = 0 \);
    \item \( \{t_0 + h_n\}_{n=1}^{\infty} \subseteq (t_0, \infty) \cap \mathbb{T} \);
    \item \( c \) is a cluster point of the sequence \( \{f(t_0 + h_n)\}_{n=1}^{\infty} \).
\end{enumerate}
The collection of all such right-sided cluster points of \( f \) at \( t_0 \) is denoted by \( \mathcal{C}_{R(t_0)}(f) \).

Similarly, the \emph{left-sided cluster points} of \( f \) at \( t_0 \) are defined by considering sequences \( \{h_n\}_{n=1} ^\infty\) such that \( t_0 - h_n \in (-\infty, t_0) \cap \mathbb{T} \), and the corresponding set is denoted by \( \mathcal{C}_{L(t_0)}(f) \).
\end{definition}
Let us now introduce a time scale associated with a fixed point \( t_0 \in \mathbb{T} \). Define
\[
\tilde{\mathbb{T}} = \{ h \in \mathbb{R} : t_0 + h \in \mathbb{T} \}.
\]

\begin{definition}
Let \( f : \mathbb{T} \to \mathbb{R} \) be a real-valued function, and let \( t_0 \in \mathbb{T}^{\kappa} \) be a point. The slope function of \( f \) at  \( t_0 \) is defined as the mapping \( \varphi_f : \tilde{\mathbb{T}} \to \mathbb{R} \) given by
\[
\varphi_f(h) = \frac{f(t_0 + h) - f(\sigma(t_0))}{h - \mu(t_0)}.
\]

%\( \varphi_f : \tilde{\mathbb{T}} \setminus \{0\} \to \mathbb{R} \) defined by
%\[
%\varphi_f(h) = \frac{f(t_0 + h) - f(\sigma(t_0))}{h-\mu(t_0)}.
%\]
\end{definition}
\begin{remark}
Although the set $\tilde{\mathbb{T}}$ and the slope function $\varphi_f$ is indeed dependent on the point $t_0$, a more precise notation would ideally reflect this relationship. Nevertheless, for the sake of simplicity and notational efficiency, we retain the current symbol, as our analysis consistently considers a fixed reference point.
\end{remark}

%\begin{lemma}
%Let $f : [t_0, t_0 + \delta)\cap{\mathbb{T}} \to \mathbb{R}$ be a real valued function, where $\delta > 0$ and $t_0 \in \mathbb{T}$ is a right-dense point. If $\mathcal{C}_{R(0)}(\varphi_f) \subset \mathbb{R}$, then $f$ is right continuous at $t_0$.
%\end{lemma}
%\begin{proof}
%Assume, for the sake of contradiction, that the function   $f$ is not right continuous at the point $t_0$.  Then there exists some $\varepsilon > 0$ such that for every positive integer $n \in \mathbb{N}$,  there exists a point $h_n\in\tilde{\mathbb{T}}$ with $0 < h_n < 1/n$  such that
%\begin{equation*}
%\left|f(t_0 + h_n) - f(t_0)\right| \geq \varepsilon.
%\end{equation*}
%Using this inequality, we see that
%\[
%\left|\varphi_f(h_n)\right| = \left| \frac{f(t_0 + h_n) - f(t_0)}{h_n} \right| \geq  \frac{\varepsilon}{h_n} > n\varepsilon.
%\]
%This inequality reveals that as  $n$ becomes arbitrarily large, the slope grows without bound. Therefore, we conclude
%\[
%\lim_{n \to \infty} \left|\varphi_f(h_n)\right| = +\infty,
%\]
%which implies that at least one of $-\infty$ and $+\infty$ is a right-sided cluster point of  $\varphi_f$  at $0$. However, this contradicts the original assumption that all right-sided cluster points of $\varphi_f$ at 0 are real numbers. Consequently, our initial assumption must be false. We conclude that the function  $f$ must indeed be right continuous at 
 %$t_0.$
%\end{proof}
\begin{definition}\label{rightcomp}
Let $f_1, f_2:(t_0,t_0+\delta)\cap{\mathbb{T}} \to \mathbb{R}$ be two real valued functions, where $\delta > 0$ and $t_0 \in \mathbb{T}$ is a right-dense point. We say that $f_1$ and $f_2$ are \textit{right complementary} at $t_0$ if they satisfy the following two conditions:
\begin{itemize}
    \item[(i)] $
    \mathcal{C}_{R(t_0)}(f_1) = \mathcal{C}_{R(t_0)}(f_2) = \{\underline{a},\overline{a}\},\,\text{where}\,\, \underline{a},\overline{a}\in \mathbb{R},\,\underline{a} <\overline{a}.
    $
    \item[(ii)] The following limits hold:
    \[
    \lim_{\substack{h \to 0^+\\ h\in\tilde{\mathbb{T}}}}\min\{f_1(t_0 + h), f_2(t_0 + h)\} =\underline{a},
    \]
    \[
     \lim_{\substack{h \to 0^+\\ h\in\tilde{\mathbb{T}}}} \max\{f_1(t_0 + h), f_2(t_0 + h)\} =\overline{a}.\]
    \end{itemize}

Similarly, the functions \( f_1 \) and \( f_2 \) are said to be \emph{left complementary} at \(t_0 \) if they fulfill the analogous criteria at \( t_0 \).
\end{definition}
\begin{definition}\label{defder}\cite{Mert1}
Let \( f: \mathbb{T} \to \mathbb{R}_{\cal{F}} \) be a fuzzy function, and let \( t_0 \in \mathbb{T}^{\kappa} \). The generalized Hukuhara \(\Delta\)-derivative of \( f \) at \( t_0 \), denoted by \( \Delta_{gH} f(t_0) \), is defined as the element of \( \mathbb{R}_{\cal{F}} \) (if it exists) with the property that for any \( \varepsilon > 0 \), there exists a neighborhood \( U_{\mathbb{T}}=(t_0- \delta, t_0 + \delta) \cap \mathbb{T}\) 
of \( t_0\)  for some \( \delta > 0 \) such that the following inequalities hold:
\begin{equation*}
D\left(f(t_0+h) \ominus_{gH} f(\sigma(t_0)), \Delta_{gH} f(t_0)(h - \mu(t_0))\right) \leq \varepsilon |h - \mu(t_0)|,
\end{equation*}
\begin{equation*}
D\left(f(\sigma(t_0)) \ominus_{gH} f(t_0-h), \Delta_{gH} f(t_0)(h + \mu(t_0))\right) \leq \varepsilon |h + \mu(t_0)|
\end{equation*}
for all \( t_0 - h, t_0 + h \in U_{\mathbb{T}} \) with \( 0 \leq h < \delta \).

It is said that \( f \) is \(\Delta_{gH}\)-differentiable at \( t_0 \) if the \(\Delta_{gH}\)-derivative of \( f \) exists at \( t_0 \). Furthermore, \( f \) is called  \(\Delta_{gH}\)-differentiable on \( \mathbb{T}^{\kappa} \) if the \(\Delta_{gH}\)-derivative of $f$ exists at every point \( t_0 \in \mathbb{T}^{\kappa} \). In this context, the function \( \Delta_{gH} f: \mathbb{T}^{\kappa} \to \mathbb{R}_{\cal{F}} \) is referred to as the \(\Delta_{gH}\)-derivative of \( f \) over the time scale \( \mathbb{T}^{\kappa} \).
\end{definition}
The generalized Hukuhara 
$\Delta$-derivative can alternatively be defined through a limit-based formulation, as described below.
\begin{definition}\label{defderlim}\cite{Mert1}
A fuzzy function \( f: \mathbb{T} \to \mathbb{R}_{\cal{F}} \) is said to be \(\Delta_{gH}\)-differentiable at \( t_0 \in \mathbb{T}^{\kappa} \) if there exists an element \( \Delta_{gH} f(t_0) \in \mathbb{R}_{\cal{F}} \) such that
\[
\Delta_{gH} f(t_0) = \lim_{\substack{h \to 0\\ h\in\tilde{\mathbb{T}}}}\frac{1}{h - \mu(t_0)} \left[ f(t_0+h)\ominus_{gH} f(\sigma(t
_0))\right].
\]
\end{definition}

As introduced in \cite{Mert1}, a fuzzy number-valued function \( f: \mathbb{T} \to \mathbb{R}_{\mathcal{F}} \) can be represented, for each \( \alpha \in [0,1] \), by an interval-valued function \( f_{\alpha}(t) = [f(t)]_{\alpha} \in \mathcal{K}_C \), where \( \mathcal{K}_C \) denotes the set of all bounded closed intervals in \( \mathbb{R} \). Each such interval  can be expressed as
\[
f_{\alpha}(t) = [f_{\alpha}^{-}(t), f_{\alpha}^{+}(t)],
\]
with \( f_{\alpha}^{-} \) and \( f_{\alpha}^{+} \) denoting the real valued lower and upper endpoint functions, respectively, defined on the time scale \( \mathbb{T} \). In that work, the authors have investigated the connection between the \(\Delta_{gH}\)-differentiability of the fuzzy function \( f \) and the \(\Delta_{gH}\)-differentiability of its associated family of interval-valued functions \(  f_{\alpha} \), and have established the following theorem.

\begin{theorem}\label{Thmlevel}\cite{Mert1}
If \( f: \mathbb{T} \to \mathbb{R}_{\cal{F}} \) is \(\Delta_{gH}\)-differentiable at \( t_0\in \mathbb{T}^{\kappa} \), then the family of interval-valued functions \( f_{\alpha}: \mathbb{T} \to \mathcal{K}_C \) is \(\Delta_{gH}\)-differentiable at \( t_0 \) uniformly in \( \alpha \in [0,1] \). Furthermore, the \(\Delta_{gH}\)-derivative of \( f_{\alpha} \) at \( t_0 \) is given by
\[
\Delta_{gH} f_{\alpha}(t_0) = [\Delta_{gH} f(t
_0)]_{\alpha}
\]
for all \( \alpha \in [0,1] \).
\end{theorem}

It can be readily shown that the converse of the above theorem also holds.

\begin{theorem}\label{reverse}
Let \( f: \mathbb{T} \to \mathbb{R}_{\mathcal{F}} \) be a fuzzy number-valued function, and let \( t_0 \in \mathbb{T}^{\kappa} \) be a point. Suppose that, for every \( \alpha \in [0,1] \),  the associated interval-valued functions \( f_{\alpha}(t) = [f(t)]_{\alpha} \) are \(\Delta_{gH}\)-differentiable at \( t_0 \), uniformly in \( \alpha \). Furthermore, assume that the generalized Hukuhara difference \( f(t_0 + h) \ominus_{gH} f(\sigma(t_0)) \) exists for all \( h \in \tilde{\mathbb{T}} \) such that \( |h| < \delta \), for some sufficiently small \( \delta > 0 \). Then the function  $f$ is \(\Delta_{gH}\)-differentiable at \( t_0 \).
\end{theorem}

\begin{remark}\label{Rmk}
The conclusions of Theorems~\ref{Thmlevel} and~\ref{reverse} remain valid for both the right-hand and left-hand \(\Delta_{gH}\)-derivatives.
\end{remark}

The following results extend Theorem~\ref{nabladerrule}, which concerns the \(\Delta\)-derivative of real-valued functions on time scales, to the more general setting of the \(\Delta_{gH}\)-derivative for fuzzy number-valued functions defined on a time scale.

\begin{theorem}\label{derivative}(\cite{Mert1})
Let \( f: \mathbb{T} \to \mathbb{R}_{\cal{F}} \) be a fuzzy number-valued function, and let \( t_0 \in \mathbb{T}^{\kappa} \) be a point. Then the following results hold:
\begin{enumerate}[label=(\roman*)]
\item If the \(\Delta_{gH}\)-derivative of \( f \) exists at \( t_0 \), then it is unique.
\item If $f$ is \(\Delta_{gH}\)-differentiable at $t_0,$ then $f$ is continuous at $t_0.$
\item If \( f \) is continuous at \( t_0 \) and \( t_0 \) is right-scattered, then \( f \) is \(\Delta_{gH}\)-differentiable at \( t_0 \), and the derivative is given by
\[
\Delta_{gH} f(t_0) = \frac{f(\sigma(t_0)) \ominus_{gH} f(t_0)}{\mu(t_0)}.
\]
\item If \( t_0 \) is right-dense, then \( f \) is \(\Delta_{gH}\)-differentiable at \( t_0 \) if and only if the following two limits exist:
\[
\lim_{\substack{h \to 0^+\\ h\in\tilde{\mathbb{T}}}} \frac{f(t_0+h) \ominus_{gH} f(t_0)}{h} \quad \text{and} \quad \lim_{\substack{h \to 0^+\\ h\in\tilde{\mathbb{T}}}}\frac{f(t_0) \ominus_{gH} f(t_0-h)}{h},
\]
and the equality
\[
\lim_{\substack{h \to 0^+\\ h\in\tilde{\mathbb{T}}}} \frac{f(t_0+h) \ominus_{gH} f(t_0)}{h} = \lim_{\substack{h \to 0^+\\ h\in\tilde{\mathbb{T}}}} \frac{f(t_0) \ominus_{gH} f(t_0-h)}{h} = \Delta_{gH} f(t_0)
\]
holds.
\end{enumerate}
\end{theorem}
\begin{theorem}\cite{Mert1}\label{rhothm}
Let \( f : \mathbb{T} \to \mathbb{R}_{\cal{F}} \) be a fuzzy number-valued function, and let \( t_0 \in \mathbb{T}^{\kappa} \) be a point. If \( f \) is \( \Delta_{gH} \)-differentiable at \( t_0 \), then the following identities hold:
\[
f(\sigma(t_0)) = f(t_0) \oplus \mu(t_0) \Delta_{gH} f(t_0),
\]
or
\[
f(t_0) = f(\sigma(t_0)) \oplus (-1) \mu(t_0) \Delta_{gH} f(t_0).
\]
\end{theorem}
\section{Main Results}
In this section, we provide a detailed analysis of the generalized Hukuhara 
\(\Delta\)-differentiability for fuzzy number-valued functions defined on time scales, expressed in terms of their associated lower and upper endpoint representations.

Our analysis is restricted to the left \(\Delta\)-derivative at points that are simultaneously right-scattered and left-dense, and the right \(\Delta\)-derivative at points that are both left-scattered and right-dense, within the time scale \(\mathbb{T}^{\kappa}\).

For the sake of notational simplicity, we introduce the following definitions, which will be used throughout the paper as needed. For each \(\alpha \in [0,1]\),  we define
\[
\lim_{\substack{h \to 0^+\\ h \in \tilde{\mathbb{T}}}} \min\left\{ \varphi_{f_{\alpha}^{-}}(h), \varphi_{f_{\alpha}^{+}}(h) \right\}:=a_{\alpha}:=\lim_{\substack{h \to 0^-\\ h \in \tilde{\mathbb{T}}}} \min\left\{ \varphi_{f_{\alpha}^{-}}(h), \varphi_{f_{\alpha}^{+}}(h) \right\},\]
\[
\lim_{\substack{h \to 0^+\\ h \in \tilde{\mathbb{T}}}} \max\left\{ \varphi_{f_{\alpha}^{-}}(h), \varphi_{f_{\alpha}^{+}}(h) \right\}:=b_{\alpha}:=\lim_{\substack{h \to 0^-\\ h \in \tilde{\mathbb{T}}}} \max\left\{ \varphi_{f_{\alpha}^{-}}(h), \varphi_{f_{\alpha}^{+}}(h) \right\} .
\]

These expressions, denoted by \(a_{\alpha}\) and \(b_{\alpha}\), will appear frequently in the analysis that follows. We encourage the reader to keep this notation in mind, as it contributes to the clarity and consistency of the exposition.

%We will refer to these expressions as \(a_{\alpha}\) and \(b_{\alpha}\), respectively, in several parts of the paper. The reader is encouraged to keep this shorthand in mind, as it will be used repeatedly for clarity and conciseness.
\begin{theorem}\label{right1}
Let \( f:\mathbb{T} \to \mathbb{R}_{\mathcal{F}} \) be a fuzzy number-valued function, and let \( t_0 \in \mathbb{T}^{\kappa} \) be a right-dense point. Suppose that \( f \) is right \(\Delta_{gH}\)-differentiable at \( t_0 \). Then, one of the following two cases must hold:

\begin{enumerate}[label=(\roman*)]
\item The right delta derivatives \( \Delta_{+} f_{\alpha}^{-}(t_0) \) and \( \Delta_{+} f_{\alpha}^{+}(t_0) \)  exist uniformly with respect to  \( \alpha \in [0,1] \). Hence, the right \(\Delta_{gH}\)-derivative of \( f \) at \( t_0 \) takes the form
\[
[\Delta_{gH_+} f(t_0)]_{\alpha} = \left[ \Delta_{+} f_{\alpha}^{-}(t_0), \Delta_{+} f_{\alpha}^{+}(t_0) \right] \]\quad \text{or} \quad
\[[\Delta_{gH_+} f(t_0)]_{\alpha} = \left[ \Delta_{+} f_{\alpha}^{+}(t_0), \Delta_{+} f_{\alpha}^{-}(t_0) \right]
\]
for all $\alpha.$%\( \alpha \in [0,1] \).

\item There exists a threshold \( \alpha_0 \in (0,1] \) such that, for all \( \alpha \in [\alpha_0,1] \), the right delta derivatives \( \Delta_{+} f_{\alpha}^{-}(t_0) \) and \( \Delta_{+} f_{\alpha}^{+}(t_0) \) exist and are equal, that is,
\[
\Delta_{+} f_{\alpha}^{-}(t_0) = \Delta_{+} f_{\alpha}^{+}(t_0) = a_{\alpha} = b_{\alpha},
\]
and this existence is uniform with respect to \( \alpha \). 

\noindent On the other hand, for \( \alpha \in [0,\alpha_0) \), the functions \( f_{\alpha}^{-} \) and \( f_{\alpha}^{+} \) are not right delta differentiable at \( t_0 \), but their associated slope functions \( \varphi_{f_{\alpha}^{-}} \) and \( \varphi_{f_{\alpha}^{+}} \) are right complementary at zero. That is, there exist constants \( a_{\alpha} < b_{\alpha} \) such that
\[
C_{R(0)}(\varphi_{f_{\alpha}^{-}}) = C_{R(0)}(\varphi_{f_{\alpha}^{+}}) = \{a_{\alpha}, b_{\alpha}\},
\]
and the limits 
\[
\lim_{\substack{h \to 0^+ \\ h \in \tilde{\mathbb{T}}}} \min\left\{ \varphi_{f_{\alpha}^{-}}(h), \varphi_{f_{\alpha}^{+}}(h) \right\} = a_{\alpha}, \quad
\lim_{\substack{h \to 0^+ \\ h \in \tilde{\mathbb{T}}}} \max\left\{ \varphi_{f_{\alpha}^{-}}(h), \varphi_{f_{\alpha}^{+}}(h) \right\} = b_{\alpha}\]
exist uniformly with respect to  $\alpha$. 

Consequently, the right \(\Delta_{gH}\)-derivative of \( f \) at \( t_0 \) is given by
\[
[\Delta_{gH_+} f(t_0)]_{\alpha} = [a_{\alpha}, b_{\alpha}]
\]
for all \( \alpha \in [0,1] \).
\end{enumerate}
\end{theorem}
\begin{proof}
Suppose that the function \( f \) is right \(\Delta_{gH}\)-differentiable at a right-dense point $t_0\in\mathbb{T}^{\kappa}.$ We consider two possible cases.

\textbf{Case 1:} The value of \( \Delta_{gH_{+}} f(t_0) \) is a crisp (real)  number.

For any \( h > 0 \) such that \( h \in \tilde{\mathbb{T}}\),  the following inequalities hold naturally: 
\begin{equation}\label{E1S}
\min\left\{ \varphi_{f_{\alpha }^{-}}(h), \varphi_{f_{\alpha }^{+}}(h) \right\}\leq \varphi_{f_{\alpha }^{-}}(h)\leq \max\left\{\varphi_{f_{\alpha }^{-}}(h), \varphi_{f_{\alpha }^{+}}(h) \right\},
\end{equation}
\begin{equation}\label{E2S}
\min\left\{ \varphi_{f_{\alpha }^{-}}(h), \varphi_{f_{\alpha }^{+}}(h) \right\}\leq \varphi_{f_{\alpha }^{+}}(h)\leq \max\left\{\varphi_{f_{\alpha }^{-}}(h), \varphi_{f_{\alpha }^{+}}(h) \right\}.
\end{equation}
Since  the function \(f \) is   assumed to be right \(\Delta_{gH}\)-differentiable at  the point \( t_0 \), it follows from  Definition \ref{defderlim} and Theorem \ref{prop2} %Definition \ref{defderlim}
that, for each $\alpha\in[0,1],$ 
the level set of the right $\Delta_{gH}$ derivative is given by
\begin{align*}
[\Delta_{gH_+} f(t_0)]_{\alpha}
&= \lim_{\substack{h \to 0^+\\ h\in\tilde{\mathbb{T}}}} \frac{1}{h} \left( [f(t_0 + h)]_{\alpha} \ominus_{gH} [f(t_0)]_{\alpha} \right) \notag \\
&= \lim_{\substack{h \to 0^+\\ h\in\tilde{\mathbb{T}}}} \frac{1}{h} \left[
\min\left\{ f_{\alpha}^{-}(t_0 + h) - f_{\alpha}^{-}(t_0),\,
           f_{\alpha}^{+}(t_0 + h) - f_{\alpha}^{+}(t_0) \right\}, \right. \notag \\
&\qquad\qquad\qquad\left.
\max\left\{ f_{\alpha}^{-}(t_0 + h) - f_{\alpha}^{-}(t_0),\,
           f_{\alpha}^{+}(t_0 + h) - f_{\alpha}^{+}(t_0) \right\} \right] \notag \\
&= \lim_{\substack{h \to 0^+\\ h\in\tilde{\mathbb{T}}}} \left[
\min\left\{ \frac{f_{\alpha}^{-}(t_0 + h) - f_{\alpha}^{-}(t_0)}{h},\,
           \frac{f_{\alpha}^{+}(t_0 + h) - f_{\alpha}^{+}(t_0)}{h} \right\}, \right. \notag \\
&\qquad\qquad\qquad\left.
\max\left\{ \frac{f_{\alpha}^{-}(t_0 + h) - f_{\alpha}^{-}(t_0)}{h},\,
           \frac{f_{\alpha}^{+}(t_0 + h) - f_{\alpha}^{+}(t_0)}{h} \right\} \right] \notag \\
    &= \left[
\lim_{\substack{h \to 0^+\\ h\in\tilde{\mathbb{T}}}}\min\left\{ \varphi_{f_{\alpha }^{-}}(h), \varphi_{f_{\alpha }^{+}}(h) \right\},\,
\lim_{\substack{h \to 0^+\\ h\in\tilde{\mathbb{T}}}} \max\left\{ \varphi_{f_{\alpha }^{-}}(h), \varphi_{f_{\alpha }^{+}}(h) \right\}
\right].
%&=[a_{\alpha}, b_{\alpha}].\label{EEE}
\end{align*}
The fact that \( \Delta_{gH_{+}} f(t_0) \) is a crisp number implies that the left and right endpoints of the \(\alpha\)-level set interval are identical. Hence, for each \( \alpha \in [0,1] \), we have
%\[[a_{\alpha}, b_{\alpha}]=
%\Delta_{gH_{+}} f(t_0).\]
\[\lim_{\substack{h \to 0^+ \\ h \in \tilde{\mathbb{T}}}} \min\left\{ \varphi_{f_\alpha^{-}}(h), \varphi_{f_\alpha^{+}}(h) \right\}
=
\lim_{\substack{h \to 0^+ \\ h \in \tilde{\mathbb{T}}}} \max\left\{ \varphi_{f_\alpha^{-}}(h), \varphi_{f_\alpha^{+}}(h) \right\}
=
\Delta_{gH_{+}} f(t_0).
\]
%Now, since \( \Delta_{gH_{+}} f(t_0) \) is a crisp number,the lower and upper bounds of the interval described above must be equal. That is, for all $\alpha\in[0,1],$
%\begin{equation}\label{EEEEE}
% \lim_{\substack{h \to 0^+\\ h\in\tilde{\mathbb{T}}}}\min\left\{ \varphi_{f_{\alpha }^{-}}(h), \varphi_{f_{\alpha }^{+}}(h) \right\} =  \lim_{\substack{h \to 0^+\\ h\in\tilde{\mathbb{T}}}}\max\left\{ \varphi_{f_{\alpha }^{-}}(h), \varphi_{f_{\alpha }^{+}}(h) \right\}=\Delta_{gH_{+}} f(t_0). 
%\end{equation}
By applying the Squeeze Theorem along with inequalities \eqref{E1S} and \eqref{E2S}, we conclude that the right-hand limits of both slope functions exist and are equal. Specifically,
\[
\Delta_{+} f_{\alpha }^{-}(t_0)=\lim_{\substack{h \to 0^+\\ h\in\tilde{\mathbb{T}}}} \varphi_{f_{\alpha }^{-}}(h) = \Delta_{gH_{+}} f(t_0),
\]
\[
\Delta_{+} f_{\alpha }^{+}(t_0)=\lim_{\substack{h \to 0^+\\ h\in\tilde{\mathbb{T}}}} \varphi_{f_{\alpha }^{+}}(h) =\Delta_{gH_{+}} f(t_0),
\]
for all $\alpha\in[0,1].$ 

\noindent Combining this result  with the definition of the Hausdorff metric, it follows that both \( f_{\alpha}^{-} \) and \( f_{\alpha}^{+} \) are right delta differentiable at \( t_0 \) for every \( \alpha \in [0,1] \), and this differentiability is uniform in \( \alpha \). Consequently, the level set  of the right  \(\Delta_{gH}\)-derivative is given by
 \[
    [\Delta_{gH_+} f(t_0)]_{\alpha} = \Delta_{gH_{+}} f(t_0)
    %\left[ \min\left\{\Delta_{+} f_{\alpha }^{-}(t_0), \Delta_{+} f_{\alpha }^{+}(t_0) \right\}, \max\left\{ \Delta_{+} f_{\alpha }^{-}(t_0), \Delta_{+} f_{\alpha }^{+}(t_0) \right\} \right]
\]
for all \( \alpha \in [0,1] \).

\textbf{Case 2:}  The value of \( \Delta_{gH_{+}} f(t_0) \) is not a crisp  number.

To analyze this situation, we consider three subcases based on the behavior of the auxiliary functions \( \varphi_{f_{\alpha}^{-}}(h) \) and \( \varphi_{f_{\alpha}^{+}}(h) \).

\textbf{Subcase 1:}  There exists a constant  \( \delta > 0 \) such that for every \( h \in \tilde{\mathbb{T}} \) with \( 0 \leq h < \delta \), the following inequality holds:
\begin{equation*}\label{EQ9}
\operatorname{len}\left( [f(t_0 + h)]_{\alpha} \right) \geq \operatorname{len}\left( [f(t_0)]_{\alpha} \right), \quad \text{for all } \alpha \in [0, 1].
\end{equation*}

\textbf{Subcase 2:} 
There exists  a constant  \( \delta> 0 \) such that for every \( h \in \tilde{\mathbb{T}} \) with \( 0 \leq h < \delta\), the following inequality  holds:
    \begin{equation*}\label{EQ10}
\operatorname{len}\left( [f(t_0 + h)]_{\alpha} \right) \leq \operatorname{len}\left( [f(t_0)]_{\alpha} \right), \quad \text{for all } \alpha \in [0, 1].
\end{equation*}

 \textbf{Subcase 3:}  For every \( n \in \mathbb{N} \), there exist \( h_n, h_n^{\prime} \in \tilde{\mathbb{T}} \) such that \( 0 < h_n, h_n^{\prime} <1/n \), and the following inequalities hold:
\begin{align*}
\operatorname{len}\left( [f(t_0 + h_n)]_{\alpha} \right) &> \operatorname{len}\left( [f(t_0)]_{\alpha} \right), \\
\operatorname{len}\left( [f(t_0 + h_n^{\prime})]_{\alpha} \right) &< \operatorname{len}\left( [f(t_0)]_{\alpha} \right),
\end{align*}
for all \( \alpha \in [0, 1] \).

\textbf{Subcase 1:}  
In this case, the inequality implies that for  sufficiently small   \( h>0 \) within  $\tilde{\mathbb{T}}$, the generalized Hukuhara difference \( f(t_0 + h) \ominus_{gH} f(t_0)\) simplifies to  the classical  Hukuhara difference  \(f(t_0 + h) \ominus_H f(t_0) \). 

Accordingly, we have
\[
\min\left\{ \varphi_{f_{\alpha}^{-}}(h), \varphi_{f_{\alpha}^{+}}(h) \right\} = \varphi_{f_{\alpha}^{-}}(h), \quad 
\max\left\{ \varphi_{f_{\alpha}^{-}}(h), \varphi_{f_{\alpha}^{+}}(h)\right \} = \varphi_{f_{\alpha}^{+}}(h)
\]
for all \( \alpha \in [0,1]. \) %and   sufficiently small 
%$h>0$ in $\tilde{\mathbb{T}}.$

Since  the function \( f \) is right \(\Delta_{gH}\)-differentiable at \( t_0 \), it follows that
\[
a_{\alpha} = \lim_{\substack{h \to 0^+\\ h\in\tilde{\mathbb{T}}}}\varphi_{f_{\alpha}^{-}}(h) = \Delta_{+} f_{\alpha}^{-}(t_0),
\]
\[
b_{\alpha} = \lim_{\substack{h \to 0^+\\ h\in\tilde{\mathbb{T}}}} \varphi_{f_{\alpha}^{+}}(h) = \Delta_{+} f_{\alpha}^{+}(t_0),
\]
for every \( \alpha \in [0,1] \). Therefore, both \( f_{\alpha}^{-} \) and \( f_{\alpha}^{+} \) are right delta differentiable at \( t_0 \), uniformly in  \( \alpha \), and we obtain
\[
    [\Delta_{gH_+} f(t_0)]_{\alpha} = \left[ \Delta_{+} f_{\alpha }^{-}(t_0),  \Delta_{+} f_{\alpha }^{+}(t_0) \right]\] 
for all \( \alpha \in [0,1] \).   
%\[
%[\Delta_{gH_+} f(t_0)]_{\alpha} %= 
%\left[ \min\left\{ \Delta_{+} f_{\alpha}^{-}(t_0), \Delta_{+} f_{\alpha}^{+}(t_0) \right\}, 
%\max\left\{ \Delta_{+} f_{\alpha}^{-}(t_0), \Delta_{+} f_{\alpha}^{+}(t_0) \right\} \right].
%\]

\textbf{Subcase 2:}  
%The proof follows analogously to Subcase (I), with reversed inequalities and the same conclusion regarding differentiability.
%For Subcase (II), the proof is similar to Subcase (I).
The argument proceeds similarly to Subcase 1, with the roles of \( f_{\alpha}^{-} \) and \( f_{\alpha}^{+} \) interchanged due to the reversed inequality. As a result, we have
\[
[\Delta_{gH_+} f(t_0)]_{\alpha} = \left[ \Delta_{+} f_{\alpha}^{+}(t_0), \Delta_{+} f_{\alpha}^{-}(t_0) \right]
\]
for all \( \alpha \in [0,1] \). 

\textbf{Subcase 3:}  
In this case, as \( h \to 0^+ \) with \( h \in \tilde{\mathbb{T}} \), the length of the \(\alpha\)-level sets of \( f(t_0 + h) \) oscillates, taking values both greater and smaller than that of \( [f(t_0)]_{\alpha} \), for all \( \alpha \in [0,1] \). Nevertheless, under the assumption that \( f \) is right \(\Delta_{gH}\)-differentiable at \( t_0 \), it follows that for every \( \alpha \in [0,1] \), the limits
\[
\lim_{n \to \infty} \varphi_{f_{\alpha}^{-}}(h_n) =a_{\alpha}, \quad 
\lim_{n \to \infty} \varphi_{f_{\alpha}^{+}}(h_n) =b_{\alpha},
\]
\[
\lim_{n \to \infty} \varphi_{f_{\alpha}^{+}}(h_n') = a_{\alpha}, \quad 
\lim_{n \to \infty} \varphi_{f_{\alpha}^{-}}(h_n') = b_{\alpha}
\]
hold. These equalities show that the points $a_\alpha$ and $b_\alpha$
belong to both cluster sets \( C_{R(0)}(\varphi_{f_{\alpha}^{-}}) \) and \( C_{R(0)}(\varphi_{f_{\alpha}^{+}}) \).

Since  $a_{\alpha}\leq b_{\alpha}$, suppose that for some $\alpha_0\in(0,1],$ $a_{\alpha_0}=b_{\alpha_0}.$ Then, it follows that for every $\alpha\in[\alpha_0,1],$ the equality $a_{\alpha}=b_{\alpha}$ holds. Consequently, the right delta derivatives  \( \Delta_{+} f_{\alpha}^{-}(t_0) \) and \( \Delta_{+} f_{\alpha}^{+}(t_0) \) exist and coincide. Specifically,  we have
\[
\Delta_{+} f_{\alpha}^{-}(t_0) = \Delta_{+} f_{\alpha}^{+}(t_0) = a_{\alpha} = b_{\alpha},
\]
and the existence of these derivatives is uniform with respect to \( \alpha \).

We now aim to show that for every \( \alpha \in [0,{\alpha_0})\), the cluster sets contain no points other than \( a_\alpha \) and \( b_\alpha \). In other words,
\[
C_{R(0)}(\varphi_{f_{\alpha}^{-}}) = \{ a_\alpha, b_\alpha \}\quad{\text{and}}\quad C_{R(0)}(\varphi_{f_{\alpha}^{+}}) = \{ a_\alpha, b_\alpha \}.
\]

Assume, for the sake of contradiction, that there exists some  $\alpha^*\in[0,{\alpha_0})$, and a third cluster point \( c\in C_{R(0)}(\varphi_{f_{\alpha^*}^{-}})  \)  such that \[ a_{\alpha^*} < c < b_{\alpha^*}. \] Then, there exists a sequence \( \{ h_n^{''} \}_{n=1}^\infty \subseteq \tilde{\mathbb{T}} \),   with each  term $h_n{''}  > 0$  and $h_n{''}  \to 0^+$ as $n \to \infty$, such that
\[
\lim_{n \to \infty} \varphi_{f_{\alpha^*}^{-}}(h_n'') = c.
\]
Hence, there exists \( N_1 \in \mathbb{N} \) such that for all \( n > N_1 \),
\begin{equation}\label{L1}
\frac{a_{\alpha^*} + c}{2} < \varphi_{f_{\alpha^*}^{-}}(h_n'') < \frac{c +b _{\alpha^*}}{2}.
\end{equation}
Since
\[
\lim_{n \to \infty} \min\left\{ \varphi_{f_{\alpha^*}^{-}}(h_n''), \varphi_{f_{\alpha^*}^{+}}(h_n'') \right\} = a_{\alpha^*},
\]
there is \( N_2 \in \mathbb{N} \) such that for all \( n > N_2 \),
\begin{equation}\label{L2}
\min\left\{ \varphi_{f_{\alpha^*}^{-}}(h_n''), \varphi_{f_{\alpha^*}^{+}}(h_n'') \right\} < \frac{a_{\alpha^*} + c}{2}.
\end{equation}
Similarly, as
\[
\lim_{n \to \infty} \max\left\{ \varphi_{f_{\alpha^*}^{-}}(h_n''), \varphi_{f_{\alpha^*}^{+}}(h_n'') \right\} = b_{\alpha^*},
\]
there is \( N_3 \in \mathbb{N} \) such that for all \( n > N_3 \),
\begin{equation}\label{L3}
\max\left\{ \varphi_{f_{\alpha^*}^{-}}(h_n''), \varphi_{f_{\alpha^*}^{+}}(h_n'') \right\} > \frac{c + b_{\alpha^*}}{2}.
\end{equation}
Combining inequalities \eqref{L1} and \eqref{L2},  we find that  for all \( n > \max\{N_1, N_2\} \),
\begin{equation}\label{L4}
\varphi_{f_{\alpha^*}^{+}}(h_n'') < \frac{a_{\alpha^*} + c}{2},
\end{equation}
whereas from  \eqref{L1} and \eqref{L3},   it follows that for all \( n > \max\{N_1, N_3\} \),
\begin{equation}\label{L5}
\varphi_{f_{\alpha^*}^{+}}(h_n'') > \frac{c + b_{\alpha^*}}{2}.
\end{equation}
For all \( n > \max\{N_1, N_2, N_3\} \), the inequalities \eqref{L4} and \eqref{L5} contradict each other, which means the point \( c  \) cannot exist. A similar argument applies to the cluster set $C_{R(0)}(\varphi_{f_{\alpha}^{+}}),$   which also contains no points other than \( a_\alpha \) and \( b_\alpha \). Therefore, we conclude that for every \( \alpha \in [0,\alpha_{0}),\)  the sets of right-sided cluster points of both functions are precisely
\[
C_{R(0)}(\varphi_{f_{\alpha}^{-}}) = \{ a_\alpha, b_\alpha \}, \quad 
C_{R(0)}(\varphi_{f_{\alpha}^{+}}) = \{ a_\alpha, b_\alpha \}.
\]
%for every \( \alpha \in [0,\alpha_{0}).\)

Since \( f \) is right \(\Delta_{gH}\)-differentiable at \( t_0 \), the \(\alpha\)-level set of its right \(\Delta_{gH}\)-derivative  is expressed as
\begin{align*}
[\Delta_{gH_+} f(t_0)]_{\alpha} 
&= \left[
\lim_{\substack{h \to 0^+ \\ h \in \tilde{\mathbb{T}}}} \min \left\{ \varphi_{f_{\alpha}^{-}}(h), \varphi_{f_{\alpha}^{+}}(h) \right\},
\lim_{\substack{h \to 0^+ \\ h \in \tilde{\mathbb{T}}}} \max \left\{ \varphi_{f_{\alpha}^{-}}(h), \varphi_{f_{\alpha}^{+}}(h) \right\}
\right] \\
&= [a_\alpha, b_\alpha],
\end{align*}
which is valid for all  \( \alpha \in [0,1] \), and  in particular for  \( \alpha \in [0,\alpha_0). \) This confirms that the functions \( \varphi_{f_{\alpha}^{-}}(h) \) and \( \varphi_{f_{\alpha}^{+}}(h) \) are right complementary at \( h = 0 \),  and  that the limits of the interval endpoints vary uniformly in
 $\alpha,$ as ensured by the structure of the Hausdorff metric.
\end{proof}
\begin{theorem}\label{right22}
Let $f:\mathbb{T} \to \mathbb{R}_{\cal{F}}$ be a fuzzy number-valued function, and let $t_0 \in \mathbb{T}^{\kappa}$ be a right-dense point. Then, $f$ is right $\Delta_{gH}$-differentiable at  $t_0$ if one of the following three cases holds:
\begin{enumerate}[label=(\roman*)]
\item The right delta derivatives \(\Delta_{+} f_{\alpha}^{-}(t_0)\) and \(\Delta_{+} f_{\alpha}^{+}(t_0)\) exist uniformly with respect to \(\alpha \in [0,1]\). Furthermore, as functions of \(\alpha\), \(\Delta_{+} f_{\alpha}^{-}(t_0)\) is monotonically increasing, while \(\Delta_{+} f_{\alpha}^{+}(t_0)\) is monotonically decreasing, and they satisfy the inequality
\[
\Delta_{+} f_{1}^{-}(t_0) \leq \Delta_{+} f_{1}^{+}(t_0).
\]
Under these conditions, the right \(\Delta_{gH}\)-derivative of \(f\) at \(t_0\) is expressed as
\[
\left[\Delta_{gH_+} f(t_0)\right]_{\alpha} = \left[\Delta_{+} f_{\alpha}^{-}(t_0), \, \Delta_{+} f_{\alpha}^{+}(t_0)\right]
\]
for every \(\alpha \in [0,1]\).

\item The right delta derivatives $\Delta_{+} f_{\alpha }^{-}(t_0)$ and $\Delta_{+} f_{\alpha }^{+}(t_0)$ exist  uniformly with respect to  $\alpha \in [0,1].$  Furthermore, as functions of  $\alpha$, $\Delta_+ f_{\alpha}^{+}(t_0)$ is monotonically increasing, while $\Delta_+ f_{\alpha}^{-}(t_0)$ is monotonically decreasing, and they satisfy the inequality $$\Delta_+ f_{1}^{+}(t_0) \leq \Delta_+ f_{1}^{-}(t_0).$$ Under these conditions, the right \(\Delta_{gH}\)-derivative of \( f \) at \( t_0 \) is expressed as
\begin{equation*}
    \left[\Delta_{gH_+} f(t_0)\right]_{\alpha} = \left[\Delta _+f_{\alpha}^{+}(t_0), \Delta_+ f_{\alpha}^{-}(t_0)\right]
\end{equation*}
for every $\alpha \in [0,1]$.

\item There exists a threshold \( \alpha_0 \in (0,1] \) such that,  for all \( \alpha \in [\alpha_0,1] \), the right delta derivatives \( \Delta_{+} f_{\alpha}^{-}(t_0) \) and \( \Delta_{+} f_{\alpha}^{+}(t_0) \) exist and are equal, that is,
\[
\Delta_{+} f_{\alpha}^{-}(t_0) = \Delta_{+} f_{\alpha}^{+}(t_0) = a_{\alpha} = b_{\alpha},
\]
and this existence is uniform with respect to \( \alpha \). 

\noindent  On the other hand, for \( \alpha \in [0,\alpha_0) \), the functions \( f_{\alpha}^{-} \) and \( f_{\alpha}^{+} \) are not right delta differentiable at \( t_0 \), but their associated slope functions \( \varphi_{f_{\alpha}^{-}} \) and \( \varphi_{f_{\alpha}^{+}} \) are right complementary at zero. Specifically, there exist constants \( a_{\alpha} < b_{\alpha} \) such that
\[
C_{R(0)}(\varphi_{f_{\alpha}^{-}}) = C_{R(0)}(\varphi_{f_{\alpha}^{+}}) = \{a_{\alpha}, b_{\alpha}\},
\]
and the limits 
\[
\lim_{\substack{h \to 0^+ \\ h \in \tilde{\mathbb{T}}}} \min\left\{ \varphi_{f_{\alpha}^{-}}(h), \varphi_{f_{\alpha}^{+}}(h) \right\} = a_{\alpha}, \quad
\lim_{\substack{h \to 0^+ \\ h \in \tilde{\mathbb{T}}}} \max\left\{ \varphi_{f_{\alpha}^{-}}(h), \varphi_{f_{\alpha}^{+}}(h) \right\} = b_{\alpha}\]
exist uniformly with respect to $\alpha$.

\noindent Furthermore, the endpoints \( a_{\alpha} \) and \( b_{\alpha} \) vary monotonically with \( \alpha \): the function \( a_{\alpha} \) increases, while \( b_{\alpha} \) decreases. These endpoints also satisfy the consistency condition \( a_1 \leq b_1 \). 

Consequently, the right \(\Delta_{gH}\)-derivative of \( f \) at \( t_0 \) is given by
\[
[\Delta_{gH_+} f(t_0)]_{\alpha} = [a_{\alpha}, b_{\alpha}]
\]
for all \( \alpha \in [0,1] \).
%\item  For each \( \alpha \in [0,1] \), the functions \( f_{\alpha}^{-} \) and \( f_{\alpha}^{+} \) fail to be right delta differentiable at the point \( t_0 \). Nonetheless, their associated slope functions, \( \varphi_{f_{\alpha}^{-}} \) and \( \varphi_{f_{\alpha}^{+}} \), exhibit a form of right sided complementary behavior at \( 0 \). In particular, there exist real numbers \( a_{\alpha} < b_{\alpha} \) such that the set of right-hand cluster points of both slope functions at zero is given by:
%\[
%C_{R(0)}(\varphi_{f_{\alpha}^{-}}) = C_{R(0)}(\varphi_{f_{\alpha}^{+}}) = \{a_{\alpha}, b_{\alpha}\}.
%\]
%Moreover, the following uniform limits with respect to \( \alpha \) are satisfied:
%\[
%\lim_{h \to 0^+} \min\left\{\varphi_{f_{\alpha}^{-}}(h), \varphi_{f_{\alpha}^{+}}(h)\right\} = a_{\alpha}, \quad
%\lim_{h \to 0^+} \max\left\{\varphi_{f_{\alpha}^{-}}(h), \varphi_{f_{\alpha}^{+}}(h)\right\} = b_{\alpha}.
%\]

%The bounds \( a_{\alpha} \) and \( b_{\alpha} \) depend monotonically on the level \( \alpha \): \( a_{\alpha} \) increases monotonically with \( \alpha \), while \( b_{\alpha} \) decreases monotonically. These bounds also satisfy the consistency condition \( a_1 \leq b_1 \).

%Consequently, the generalized Hukuhara right delta derivative of \( f \) at \( t_0 \), in terms of its level set representation, is given by:
%\[
%\left[\Delta_{gH_+}f(t_0)\right]_{\alpha} = [a_{\alpha}, b_{\alpha}],
%\]
%for all \( \alpha \in [0,1] \).
\end{enumerate}
\end{theorem}
\begin{proof}
Assume that Case (i) holds. Since the right delta derivatives \( \Delta_+ f_\alpha^-(t_0) \) and \( \Delta_+ f_\alpha^+(t_0) \) exist  uniformly in  \( \alpha \in [0,1] \), it follows from the definition of the right delta derivative that  for all sufficiently small  \( h > 0 \) with \( h \in \tilde{\mathbb{T}} \),  the following asymptotic expansions hold uniformly with respect to \( \alpha \):
\[
f_\alpha^-(t_0 + h) - f_\alpha^-(t_0) = h \Delta_+ f_\alpha^-(t_0) + o(h),
\]
\[
f_\alpha^+(t_0 + h) - f_\alpha^+(t_0) = h \Delta_+ f_\alpha^+(t_0) + o(h).
\]
Subtracting these two expressions yields
\[
f_\alpha^-(t_0 + h) - f_\alpha^-(t_0) - \left( f_\alpha^+(t_0 + h) - f_\alpha^+(t_0) \right)
= hA(\alpha) + o(h),
\]
where we define 
\[
A(\alpha) := \Delta_+ f_\alpha^-(t_0) - \Delta_+ f_\alpha^+(t_0).
\]

Given that \( \Delta_+ f_{\alpha}^{-}(t_0) \) is monotonically increasing and \( \Delta_+ f_{\alpha}^{+}(t_0) \) is monotonically decreasing with respect to \( \alpha \in [0,1] \), and that
\[
\Delta_+ f_{1}^{-}(t_0) \leq \Delta_+ f_{1}^{+}(t_0),
\]
it follows that
\[
\Delta_+ f_{\alpha}^{-}(t_0) \leq \Delta_+ f_{\alpha}^{+}(t_0) \quad \text{for all } \alpha \in [0,1],
\]
so $A(\alpha)\leq 0$ for all $\alpha \in [0,1].$

We now distinguish two subcases.

{\bf{Case 1:}} \( A(\alpha)< 0 \) for some $\alpha\in[0,1].$

Since \( h > 0 \) and \( A(\alpha) < 0 \), it follows that \( hA(\alpha) < 0 \). 
Also, because \(\frac{ o(h)}{h} \to 0 \) as \( h \to 0^+ \) within $\tilde{\mathbb{T}}$, there exists \( \delta > 0 \) such that for all $h\in\tilde{\mathbb{T}}$ with \( 0 < h < \delta \),  we have
\[
|o(h)| < \frac{|A(\alpha)|}{2} h.
\]
Thus,
\[
hA(\alpha) + o(h) < hA(\alpha) + \frac{|A(\alpha)|}{2} h =  \frac{A(\alpha)}{2}h < 0.
\]
This leads to the inequality
\[
f_\alpha^-(t_0 + h) - f_\alpha^-(t_0) < f_\alpha^+(t_0 + h) - f_\alpha^+(t_0)
\]
holding for all sufficiently small \( h > 0 \) within $\tilde{\mathbb{T}}$.%, and for all $\alpha\in[0,1].$

{\bf{Case 2:}} \( A(\alpha) = 0 \) for some $\alpha\in[0,1].$

 In this situation,  the expression  simplifies to
\[
f_\alpha^-(t_0 + h) - f_\alpha^-(t_0) - \left( f_\alpha^+(t_0 + h) - f_\alpha^+(t_0) \right)
= o(h).
\]
According to the definition of little-$o$ notation, for any  \( \varepsilon > 0 \), there exists \( \delta > 0 \) such that 
for all  $h\in\tilde{\mathbb{T}}$  with \( 0 < h < \delta \), we have 
$$ |o(h)| < \varepsilon h.$$
Set $\delta_0:=\min\{1,\delta\}.$
Then for all $h\in\tilde{\mathbb{T}}$ satisfying
$0<h<\delta_0$,  it follows that
$$ |o(h)| < \varepsilon .$$
Therefore,  for such $h$, we get
\[
f_\alpha^-(t_0 + h) - f_\alpha^-(t_0)<f_\alpha^+(t_0 + h)-f_\alpha^+(t_0)+\varepsilon.
\]
Taking the limit as  \( \varepsilon \to 0^+ \), we obtain the inequality
\[
f_\alpha^-(t_0 + h) - f_\alpha^-(t_0) \leq f_\alpha^+(t_0 + h) - f_\alpha^+(t_0),
\]
valid for    %\( \alpha \in [0,1] \) and 
all $0<h<\delta_0$ within $\tilde{\mathbb{T}}$.

\noindent In summary, regardless of whether \( A(\alpha) < 0 \) or \( A (\alpha)= 0 \),    the inequality
\[
f_\alpha^-(t_0 + h) - f_\alpha^-(t_0) \leq f_\alpha^+(t_0 + h) - f_\alpha^+(t_0)
\]
holds for all \( \alpha \in [0,1] \),  provided  \( h > 0 \) is  sufficiently small in $\tilde{\mathbb{T}}$. This ensures that the Hukuhara difference of the $\alpha$-level sets,
\[
[f(t_0 + h)]_\alpha \ominus_H [f(t_0)]_\alpha
\]
defines a valid interval for each  \( \alpha\in[0,1] \).

Now, let \( \varepsilon > 0 \) be arbitrary.  Since the right delta derivatives exist uniformly with respect to $\alpha$, there exists \( \delta > 0 \) such that for all \( h \in\tilde{\mathbb{T}}\) with \( 0 < h < \delta \), and all \( \alpha \in [0,1] \), we have
\[
\left| \frac{f_\alpha^{-}(t_0 + h) - f_\alpha^{-}(t_0)}{h} - \Delta_+ f_\alpha^{-}(t_0) \right| < \varepsilon.
\]
Previously, we considered the slope  function  $\varphi_{f_\alpha^-}(h)$
as a function of \( h \). To analyze  how this slope  varies with respect to  the parameter \( \alpha \), we now fix \(0<h<\delta \) in $\tilde{\mathbb T}$ and define
\[
\Phi_h(\alpha) := \frac{f_\alpha^-(t_0 + h) - f_\alpha^-(t_0)}{h}.
\]
Let  \(\alpha_1, \alpha_2 \in [0,1]\) with \(\alpha_1 < \alpha_2\). From the earlier inequality, we deduce
\[
\Phi_h(\alpha_1) < \Delta_+ f_{\alpha_1}^{-}(t_0) + \varepsilon, \quad
\Phi_h(\alpha_2) > \Delta_+ f_{\alpha_2}^{-}(t_0) - \varepsilon.
\]
Subtracting the two gives
\[
\Phi_h(\alpha_2) - \Phi_h(\alpha_1) > \Delta_+ f_{\alpha_2}^{-}(t_0) - \Delta_+ f_{\alpha_1}^{-}(t_0) - 2\varepsilon.
\]
Since \( \Delta_+ f_\alpha^{-}(t_0)\) is monotonically increasing in $\alpha$, we know
\[
\Delta_+ f_{\alpha_2}^{-}(t_0) - \Delta_+ f_{\alpha_1}^{-}(t_0) \geq 0,
\]
and hence
\[
\Phi_h(\alpha_2)>\Phi_h(\alpha_1)-2\varepsilon.
\]
As this inequality holds for any \(\varepsilon > 0\), taking the limit as \(\varepsilon \to 0^+\) yields
\[
\Phi_h(\alpha_2)\geq \Phi_h(\alpha_1).
\]
This implies that
$$f_{\alpha_{2}}^{-}(t_0 + h) - f_{\alpha_{2}}^{-}(t_0)\geq f_{\alpha_{1}}^{-}(t_0 + h) - f_{\alpha_{1}}^{-}(t_0)$$
for all sufficiently  small $h>0$ in $\tilde{\mathbb{T}}$. That is, the map
 $$\alpha \mapsto f_\alpha^{-}(t_0 + h) - f_\alpha^{-}(t_0)$$ is monotonically increasing. Similarly, using the fact that  \( \Delta_+ f_\alpha^{+}(t_0)\) is monotonically decreasing  and applying the same reasoning, we  conclude that
$$f_\alpha^{+}(t_0 + h) - f_\alpha^{+}(t_0)$$ is monotonically decreasing in $\alpha$.

Therefore, the generalized Hukuhara difference
 $$ f(t_0 + h) \ominus_{gH} f(t_0) $$ 
 exists for all sufficiently small \( h > 0 \) in $\tilde{\mathbb{T}}$.

Since the right delta derivatives $\Delta_+ f_{\alpha}^{-}(t_0)$ and $\Delta_+ f_{\alpha}^{+}(t_0)$ exist for every $\alpha \in [0,1]$, the associated auxiliary functions $\varphi_{f_{\alpha}^{-}}$ and $\varphi_{f_{\alpha}^{+}}$ satisfy the limits
\begin{equation*}\label{limlim}
\lim_{\substack{h \to 0^+ \\ h \in \tilde{\mathbb{T}}}} \varphi_{f_{\alpha }^{-}}(h)= \Delta_{+} f_{\alpha }^{-}(t_0), \quad \text{and} \quad \lim_{\substack{h \to 0^+ \\ h \in \tilde{\mathbb{T}}}} \varphi_{f_{\alpha }^{+}}(h) = \Delta_{+} f_{\alpha }^{+}(t_0).
\end{equation*}
%As a result, the limit of the pointwise minimum of the two auxiliary functions equals the minimum of their respective limits, that is,
%\begin{align*}
%\lim_{\substack{h \to 0^+ \\ h \in \tilde{\mathbb{T}}}}\min\{ \varphi_{f_{\alpha }^{-}}(h), \varphi_{f_{\alpha }^{+}}(h) \}
%= \min\{\lim_{\substack{h \to 0^+ \\ h \in \tilde{\mathbb{T}}}} \varphi_{f_{\alpha }^{-}}(h), \lim_{\substack{h \to 0^+ \\ h \in \tilde{\mathbb{T}}}}\varphi_{f_{\alpha }^{+}}(h)\}&=\min\{\Delta_{+} f_{\alpha }^{-}(t_0),\Delta_{+} f_{\alpha }^{+}(t_0)\}\\
%&=\Delta_{+} f_{\alpha }^{-}(t_0).
%\end{align*}
%and similarly,
%Similarly, the limit of the pointwise maximum is given by
%\begin{align*}
%\lim_{\substack{h \to 0^+ \\ h \in \tilde{\mathbb{T}}}}\max\{ \varphi_{f_{\alpha }^{-}}(h), \varphi_{f_{\alpha }^{+}}(h) \}
%= \max\{\lim_{\substack{h \to 0^+ \\ h \in \tilde{\mathbb{T}}}} \varphi_{f_{\alpha }^{-}}(h), \lim_{\substack{h \to 0^+ \\ h \in \tilde{\mathbb{T}}}}\varphi_{f_{\alpha }^{+}}(h)\}&=\max\{\Delta_{+} f_{\alpha }^{-}(t_0),\Delta_{+} f_{\alpha }^{+}(t_0)\}\\&=\Delta_{+} f_{\alpha }^{+}(t_0).
%\end{align*}
Because $f^{-}_\alpha$ and $f^{+}_\alpha$ are right delta differentiable at $t_0,$ uniformly in $\alpha $, these limits  also hold uniformly over $\alpha\in[0,1].$

Let $\{h_n\}_{n=1}^{\infty} \subseteq \tilde{\mathbb{T}}$ be a sequence with each  $h_n > 0$  and $h_n \to 0^+$ as $n \to \infty$.  Then  these limits can be expressed as
\[
\Delta_+ f^{-}_\alpha(t_0) = \lim_{n \to \infty} \varphi_{f_{\alpha }^{-}}(h_n)
, \quad \text{and} \quad 
\Delta_+ f^{+}_\alpha(t_0) = \lim_{n \to \infty} \varphi_{f_{\alpha }^{+}}(h_n).
\]
This shows that both  $\Delta_+ f^{-}_\alpha(t_0)$ and $\Delta_+ f^{+}_\alpha(t_0)$ are  uniform limits of sequences of left-continuous functions with respect to 
 $\alpha \in (0,1]$,  hence they  themselves are left-continuous on $(0,1].$
 %Consequently, they are themselves are left continuous for $\alpha \in (0,1]$.
 A similar argument gives right continuity at  $\alpha = 0$.
 %A similar argument establishes right continuity at $\alpha = 0$.

\noindent  Thus, by Theorem~\ref{Bede1}, the interval-valued function
 %By applying  Theorem~\ref{Bede1}, we conclude that the interval-valued function
\[
\left[ \Delta_+ f_{\alpha}^{-}(t_0), \Delta_+ f_{\alpha}^{+}(t_0) \right]
\]
defines a valid fuzzy number. 

It further follows that the limit 
\begin{align*}
\Delta_{gH_+} f_{\alpha}(t_0) 
&= \lim_{\substack{h \to 0^+ \\ h \in \tilde{\mathbb{T}}}} \frac{1}{h} \left( [f(t_0 + h)]_{\alpha} \ominus_{gH} [f(t_0)]_{\alpha} \right) \\
&= \lim_{\substack{h \to 0^+ \\ h \in \tilde{\mathbb{T}}}} \frac{1}{h} [ f_{\alpha}^{-}(t_0 + h) - f_{\alpha}^{-}(t_0),\,
           f_{\alpha}^{+}(t_0 + h) - f_{\alpha}^{+}(t_0)]\\
&=\lim_{\substack{h \to 0^+ \\ h \in \tilde{\mathbb{T}}}}\left[ \frac{f_{\alpha}^{-}(t_0 + h) - f_{\alpha}^{-}(t_0)}{h},\,
           \frac{f_{\alpha}^{+}(t_0 + h) - f_{\alpha}^{+}(t_0)}{h}\right]\\
%&= \lim_{\substack{h \to 0^+ \\ h \in \tilde{\mathbb{T}}}} \frac{1}{h} \left[
%\min\left\{ f_{\alpha}^{-}(t_0 + h) - f_{\alpha}^{-}(t_0),\,
%           f_{\alpha}^{+}(t_0 + h) - f_{\alpha}^{+}(t_0) \right\}, \right. \\
%&\qquad\qquad\qquad\left.
%\max\left\{ f_{\alpha}^{-}(t_0 + h) - f_{\alpha}^{-}(t_0),\,
  %         f_{\alpha}^{+}(t_0 + h) - f_{\alpha}^{+}(t_0) \right\} \right] \\
%&= \lim_{\substack{h \to 0^+ \\ h \in \tilde{\mathbb{T}}}} \left[
%\min\left\{ \frac{f_{\alpha}^{-}(t_0 + h) - f_{\alpha}^{-}(t_0)}{h},\,
%           \frac{f_{\alpha}^{+}(t_0 + h) - f_{\alpha}^{+}(t_0)}{h} \right\}, \right. \\
%&\qquad\qquad\qquad\left.
%\max\left\{ \frac{f_{\alpha}^{-}(t_0 + h) - f_{\alpha}^{-}(t_0)}{h},\,
%           \frac{f_{\alpha}^{+}(t_0 + h) - f_{\alpha}^{+}(t_0)}{h} \right\} \right] \\
    &=\left[
\lim_{\substack{h \to 0^+ \\ h \in \tilde{\mathbb{T}}}} \varphi_{f_{\alpha }^{-}}(h),\,
\lim_{\substack{h \to 0^+ \\ h \in \tilde{\mathbb{T}}}}   \varphi_{f_{\alpha }^{+}}(h) 
\right] \\
&= \left[\Delta_{+} f_{\alpha }^{-}(t_0), \Delta_{+} f_{\alpha }^{+}(t_0)\right]
\end{align*}
exists for all \(\alpha \in [0,1]\), and it exists uniformly in \(\alpha\).
 Hence, by Remark  \ref{Rmk}, the function $f$ is right 
$\Delta_{gH}$-differentiable at $t_0,$  and
\begin{equation*}
    \left[\Delta_{{gH}_{+}} f(t_0)\right]_{\alpha} = \left[\Delta_+ f_{\alpha}^{-}(t_0), \Delta_+ f_{\alpha}^{+}(t_0)\right]
\end{equation*}
for all $\alpha \in [0,1]$. The reasoning for Case (ii) proceeds in a similar manner.

Now consider Case (iii). For \( \alpha \in [\alpha_0, 1] \),  by the reasoning used in  Case (i), the generalized Hukuhara difference
    \[
    f(t_0 + h) \ominus_{gH} f(t_0)
    \]
    exists for all sufficiently small \( h > 0 \) in $\tilde{\mathbb{T}}$. On the other hand, for \( \alpha \in [0, \alpha_0) \), although the endpoint functions \( f_\alpha^- \) and \( f_\alpha^+ \) are not right delta differentiable at $t_0$, the uniform existence with respect to $\alpha$ of the limits of the minimum and maximum of the slope functions ensures that the generalized Hukuhara difference at level $\alpha,$
    \[[f(t_0 + h)]_\alpha \ominus_{gH} [f(t_0)]_\alpha\]
    forms a valid interval for sufficiently small \( h > 0 \) in $\tilde{\mathbb{T}}$. Depending on the slope behavior, this interval is either
\[
    [f_\alpha^-(t_0 + h) - f_\alpha^-(t_0),\, f_\alpha^+(t_0 + h) - f_\alpha^+(t_0)]
    \]
    or
     \[
     [f_\alpha^+(t_0 + h) - f_\alpha^+(t_0),\, f_\alpha^-(t_0 + h) - f_\alpha^-(t_0)].
    \]
%, with the limiting behavior of slope functions determining the derivative interval.
Therefore, the generalized Hukuhara difference
 $$ f(t_0 + h) \ominus_{gH} f(t_0) $$ 
 exists for all sufficiently small \( h > 0 \) in $\tilde{\mathbb{T}}$.

Assuming that the following limits exist uniformly for all $\alpha\in[0,1],$
\[\lim_{\substack{h \to 0^+ \\ h \in \tilde{\mathbb{T}}}} \min\left\{ \varphi_{f_{\alpha}^{-}}(h), \varphi_{f_{\alpha}^{+}}(h) \right\} = a_{\alpha},\quad
\lim_{\substack{h \to 0^+ \\ h \in \tilde{\mathbb{T}}}} \max\left\{ \varphi_{f_{\alpha}^{-}}(h), \varphi_{f_{\alpha}^{+}}(h) \right\} = b_{\alpha},
\]
we conclude that the function \( f_{\alpha} \) is right \(\Delta_{gH}\)-differentiable at \( t_0 \), and this differentiability holds uniformly with respect to \( \alpha \). Furthermore, because these limits exist uniformly in \( \alpha \),  an argument analogous to  Case (i) shows that the endpoint functions  \( a_{\alpha} \) and \( b_{\alpha} \) are left-continuous on  \( (0,1] \) and right-continuous at \( \alpha = 0 \). Hence, by Theorem~\ref{Bede1}, the family of intervals \( [a_{\alpha}, b_{\alpha}] \)  constitutes a valid fuzzy number.
 
Additionally, the following uniform limit exists and defines the right \(\Delta_{gH}\)-derivative of $f_{\alpha}$ at $t_0$:
\begin{align*}
\Delta_{gH_+} f_{\alpha}(t_0) 
&= \lim_{\substack{h \to 0^+ \\ h \in \tilde{\mathbb{T}}}} \frac{1}{h} \left( [f(t_0 + h)]_{\alpha} \ominus_{gH} [f(t_0)]_{\alpha} \right) \\
&= \lim_{\substack{h \to 0^+ \\ h \in \tilde{\mathbb{T}}}} \frac{1}{h} \left[
\min\left\{ f_{\alpha}^{-}(t_0 + h) - f_{\alpha}^{-}(t_0),\,
           f_{\alpha}^{+}(t_0 + h) - f_{\alpha}^{+}(t_0) \right\}, \right. \\
&\qquad\qquad\qquad\left.
\max\left\{ f_{\alpha}^{-}(t_0 + h) - f_{\alpha}^{-}(t_0),\,
           f_{\alpha}^{+}(t_0 + h) - f_{\alpha}^{+}(t_0) \right\} \right] \\
&= \lim_{\substack{h \to 0^+ \\ h \in \tilde{\mathbb{T}}}} \left[
\min\left\{ \frac{f_{\alpha}^{-}(t_0 + h) - f_{\alpha}^{-}(t_0)}{h},\,
           \frac{f_{\alpha}^{+}(t_0 + h) - f_{\alpha}^{+}(t_0)}{h} \right\}, \right. \\
&\qquad\qquad\qquad\left.
\max\left\{ \frac{f_{\alpha}^{-}(t_0 + h) - f_{\alpha}^{-}(t_0)}{h},\,
           \frac{f_{\alpha}^{+}(t_0 + h) - f_{\alpha}^{+}(t_0)}{h} \right\} \right] \\
&= \left[
\lim_{\substack{h \to 0^+ \\ h \in \tilde{\mathbb{T}}}} \min\left\{ \varphi_{f_{\alpha }^{-}}(h), \varphi_{f_{\alpha }^{+}}(h) \right\},\,
\lim_{\substack{h \to 0^+ \\ h \in \tilde{\mathbb{T}}}} \max\left\{ \varphi_{f_{\alpha }^{-}}(h), \varphi_{f_{\alpha }^{+}}(h) \right\}
\right] \\
&= [a_{\alpha}, b_{\alpha}].
\end{align*}
Thus, by Remark  \ref{Rmk}, the function $f$ is right  $\Delta_{gH}$-differentiable at $t_0,$   with its derivative given level-wise by
\[
\left[\Delta_{gH_+} f(t_0)\right]_\alpha = 
\begin{cases}
[a_\alpha, b_\alpha], & \alpha \in [0, \alpha_0), \\
[a_\alpha, a_\alpha], &\alpha \in [\alpha_0, 1],
\end{cases}
\]
where in the latter case the derivative reduces to a singleton interval.
\end{proof}

Similarly to  Theorems \ref{right1} and \ref{right22}, we can also establish the following results for a left-dense point \( t_0 \in \mathbb{T}^{\kappa} \).

\begin{theorem}\label{left1}
Let \( f:\mathbb{T} \to \mathbb{R}_{\mathcal{F}} \) be a fuzzy number-valued function, and let \( t_0 \in \mathbb{T}^{\kappa} \) be a left-dense point. Suppose that \( f \) is left \(\Delta_{gH}\)-differentiable at \( t_0 \). Then, one of the following two cases must hold:

\begin{enumerate}[label=(\roman*)]
\item The left delta derivatives \( \Delta_{-} f_{\alpha}^{-}(t_0) \) and \( \Delta_{-} f_{\alpha}^{+}(t_0) \)  exist uniformly with respect to  \( \alpha \in [0,1] \).  Hence, the left \(\Delta_{gH}\)-derivative of \( f \) at \( t_0 \) takes the form
  \[
    [\Delta_{gH_-} f(t_0)]_{\alpha} = \left[\Delta_{-} f_{\alpha }^{-}(t_0), \Delta_{-} f_{\alpha }^{+}(t_0) \right]
    \]
    or
    \[
    [\Delta_{gH_-} f(t_0)]_{\alpha} = \left[\Delta_{-} f_{\alpha }^{+}(t_0), \Delta_{-} f_{\alpha }^{-}(t_0) \right]
    \]
for all \( \alpha\).

\item There exists a threshold \( \alpha_0 \in (0,1] \) such that, for all \( \alpha \in [\alpha_0,1] \), the left delta derivatives \( \Delta_{-} f_{\alpha}^{-}(t_0) \) and \( \Delta_{-} f_{\alpha}^{+}(t_0) \) exist and are equal, that is,
\[
\Delta_{-} f_{\alpha}^{-}(t_0) = \Delta_{-} f_{\alpha}^{+}(t_0) = a_{\alpha} = b_{\alpha},
\]
and this existence is uniform with respect to  \( \alpha \).

\noindent On the other hand, for \( \alpha \in [0,\alpha_0) \), the functions \( f_{\alpha}^{-} \) and \( f_{\alpha}^{+} \) are not left delta differentiable at \( t_0 \), but their associated slope functions \( \varphi_{f_{\alpha}^{-}} \) and \( \varphi_{f_{\alpha}^{+}} \) are left complementary at zero. That is, there exist constants \( a_{\alpha} < b_{\alpha} \) such that
\[
    C_{L(0)}(\varphi_{f_{\alpha }^{-}}) = C_{L(0)}(\varphi_{f_{\alpha }^{+}}) = \{a_{\alpha}, b_{\alpha}\}\]  and the limits 
    \[\lim_{\substack{h \to 0^- \\ h \in \tilde{\mathbb{T}}}} \min\left\{\varphi_{f_{\alpha }^{-}}(h), \varphi_{f_{\alpha }^{+}}(h)\right\} =a_{\alpha},\quad 
    \lim_{\substack{h \to 0^- \\ h \in \tilde{\mathbb{T}}}} \max\left\{\varphi_{f_{\alpha }^{-}}(h), \varphi_{f_{\alpha }^{+}}(h)\right\} =b_{\alpha}
    \] 
exist uniformly with respect to $\alpha$. 

\noindent Consequently, the left \(\Delta_{gH}\)-derivative of \( f \) at \( t_0 \) is given by
    $$[\Delta_{gH_-} f(t_0)]_{\alpha }=[a_{\alpha}, b_{\alpha}]$$
for all $\alpha \in [0,1]$.
\end{enumerate}
\end{theorem}

\begin{theorem}\label{right2}
Let $f:\mathbb{T} \to \mathbb{R}_{\cal{F}}$ be a fuzzy number-valued function, and let $t_0 \in \mathbb{T}^{\kappa}$ be a left-dense point. Then, $f$ is left $\Delta_{gH}$-differentiable at  $t_0$ if one of the following three cases holds:
\begin{enumerate}[label=(\roman*)]
 \item The left delta derivatives $\Delta_{-} f_{\alpha }^{-}(t_0)$ and $\Delta_{-} f_{\alpha }^{+}(t_0)$ exist uniformly with respect to  $\alpha \in [0,1].$  Furthermore, as functions of $\alpha$, $\Delta_-f_{\alpha}^{-}(t_0)$ is monotonically increasing, while $\Delta_- f_{\alpha}^{+}(t_0)$ is monotonically decreasing, and they satisfy the inequality $$\Delta_- f_{1}^{-}(t_0) \leq \Delta_- f_{1}^{+}(t_0).$$Under these conditions, the left \(\Delta_{gH}\)-derivative of \( f \) at \( t_0 \) is expressed as
\begin{equation*}
    \left[\Delta_{gH_-} f(t_0)\right]_{\alpha} = \left[\Delta _-f_{\alpha}^{-}(t_0), \Delta_- f_{\alpha}^{+}(t_0)\right]
\end{equation*}
for all $\alpha \in [0,1].$

\item  The left delta derivatives $\Delta_{-} f_{\alpha }^{-}(t_0)$ and $\Delta_{-} f_{\alpha }^{+}(t_0)$ exist   uniformly with respect to  $\alpha \in [0,1].$  Furthermore, as functions of $\alpha$, $\Delta_- f_{\alpha}^{+}(t_0)$ is monotonically increasing, while $\Delta_- f_{\alpha}^{-}(t_0)$ is monotonically decreasing, and they satisfy the inequality $$\Delta_- f_{1}^{+}(t_0) \leq \Delta_- f_{1}^{-}(t_0).$$ 

\noindent Under these conditions, the left \(\Delta_{gH}\)-derivative of \( f \) at \( t_0 \) is expressed as
\begin{equation*}
    \left[\Delta_{gH_-} f(t_0)\right]_{\alpha} = \left[\Delta _-f_{\alpha}^{+}(t_0), \Delta_- f_{\alpha}^{-}(t_0)\right]
\end{equation*}
for all $\alpha \in [0,1]$.

\item There exists a threshold \( \alpha_0 \in (0,1] \) such that, for all \( \alpha \in [\alpha_0,1] \), the left delta derivatives \( \Delta_{-} f_{\alpha}^{-}(t_0) \) and \( \Delta_{-} f_{\alpha}^{+}(t_0) \) exist and are equal, that is,
\[
\Delta_{-} f_{\alpha}^{-}(t_0) = \Delta_{-} f_{\alpha}^{+}(t_0) = a_{\alpha} = b_{\alpha},
\]
and this existence is uniform with respect to \( \alpha \). 

\noindent On the other hand, for \( \alpha \in [0,\alpha_0) \), the functions \( f_{\alpha}^{-} \) and \( f_{\alpha}^{+} \) are not left delta differentiable at \( t_0 \), but their associated slope functions \( \varphi_{f_{\alpha}^{-}} \) and \( \varphi_{f_{\alpha}^{+}} \) are left complementary at zero. Specifically, there exist constants \( a_{\alpha} < b_{\alpha} \) such that
\[
    C_{L(0)}(\varphi_{f_{\alpha }^{-}}) = C_{L(0)}(\varphi_{f_{\alpha }^{+}}) = \{a_{\alpha}, b_{\alpha}\},\]  and the limits 
    \[\lim_{\substack{h \to 0^- \\ h \in \tilde{\mathbb{T}}}} \min\left\{\varphi_{f_{\alpha }^{-}}(h), \varphi_{f_{\alpha }^{+}}(h)\right\} =a_{\alpha},\quad 
    \lim_{\substack{h \to 0^- \\ h \in \tilde{\mathbb{T}}}} \max\left\{\varphi_{f_{\alpha }^{-}}(h), \varphi_{f_{\alpha }^{+}}(h)\right\} =b_{\alpha}
    \] 
exist uniformly with respect to $\alpha$.

\noindent Furthermore, the endpoints \( a_{\alpha} \) and \( b_{\alpha} \) vary monotonically with  \( \alpha \): the function \( a_{\alpha} \) increases, while \( b_{\alpha} \) decreases. These  endpoints also satisfy the consistency condition \( a_1 \leq b_1 \). 

\noindent Consequently, the left \(\Delta_{gH}\)-derivative of \( f \) at \( t_0 \) is given by
\[
\left[\Delta_{gH_-}f(t_0)\right]_{\alpha} = [a_{\alpha},b_{\alpha}]
\]
for all $\alpha \in [0,1]$.
\end{enumerate}
\end{theorem}

By combining the theorems stated above, we obtain the following results, which characterizes the gH-differentiability of a fuzzy number-valued function at a dense point \( t_0 \in \mathbb{T}^{\kappa} \).

\begin{theorem}\label{der1}
Let $f:\mathbb{T} \to \mathbb{R}_{\cal{F}}$ be a fuzzy number-valued function, and let $t_0 \in \mathbb{T}^{\kappa}$ be a  dense  point.   If $f$ is  $\Delta_{gH}$-differentiable at  $t_0$, then one of the following cases holds:
\begin{enumerate}[label=(\roman*)]

\item The functions \( f_{\alpha}^{-} \) and \( f_{\alpha}^{+} \) are uniformly \( \Delta \)-differentiable at \( t_0 \) with respect to \( \alpha \in [0,1] \). Then, for each \( \alpha \in [0,1] \), the \( \alpha \)-level set of the \( \Delta_{gH} \)-derivative of \( f \) at \( t_0 \) is an interval given by either
\[
[\Delta_{gH} f(t_0)]_{\alpha} = \left[\Delta f_{\alpha}^{-}(t_0), \Delta f_{\alpha}^{+}(t_0)\right]
\]
or
\[
[\Delta_{gH} f(t_0)]_{\alpha} = \left[\Delta f_{\alpha}^{+}(t_0), \Delta f_{\alpha}^{-}(t_0)\right].
\]

 \item The functions $f_{\alpha }^{-}$ and $f_{\alpha }^{+}$ are uniformly $\Delta$-differentiable at $t_0$   with respect to $\alpha \in [0,1]$. Then, for all \( \alpha \in [0,1] \), the \( \alpha \)-level set of the \( \Delta_{gH} \)-derivative of \( f \) at \( t_0 \) is an interval given by either
 \begin{equation*}
        [\Delta_{gH} f(t_0)]_{\alpha} = \left[\Delta f_{\alpha }^{-}(t_0), \Delta f_{\alpha }^{+}(t_0)\right]
    \end{equation*}
    or
    \begin{equation*}
        [\Delta_{gH} f(t_0)]_{\alpha} = \left[\Delta f_{\alpha }^{+}(t_0), \Delta f_{\alpha }^{-}(t_0).\right]
    \end{equation*}
    
\item   The one sided derivatives  $\Delta_{-} f_{\alpha }^{-}(t_0),$ $\Delta_{+} f_{\alpha }^{-}(t_0),$ $ \Delta_{-} f_{\alpha }^{+}(t_0),$ and $\Delta_{+} f_{\alpha }^{+}(t_0)$ exist uniformly with respect to $\alpha\in[0,1]$. Furthermore, the function  
$\Delta_{-} f_{\alpha }^{-}(t_0)=\Delta_{+} f_{\alpha }^{+}(t_0)$ is monotonically increasing with respect to $\alpha,$ while $\Delta_{+} f_{\alpha }^{-}(t_0)=\Delta_{-} f_{\alpha }^{+}(t_0)$ is monotonically decreasing.  In addition, the inequality $\Delta_{-} f_{1}^{-}(t_0)\leq \Delta_{-} f_{1}^{+}(t_0)$ holds.  In this setting, we have
$$\left[\Delta_{gH}f(t_0)\right]_{\alpha}=[\Delta_{-} f_{\alpha }^{-}(t_0),\Delta_{-} f_{\alpha }^{+}(t_0)]=[\Delta_{+} f_{\alpha }^{+}(t_0),\Delta_{+} f_{\alpha }^{-}(t_0)]$$
for all $\alpha\in[0,1].$

\item The one sided derivatives  $\Delta_{-} f_{\alpha}^{-}(t_0),$ $\Delta_{+} f_{\alpha}^{-}(t_0)$, $\Delta_{-} f_{\alpha}^{+}(t_0)$, and $\Delta_{+} f_{\alpha}^{+}(t_0)$ exist  uniformly with respect to $\alpha \in [0,1]$.   Furthermore, the function  $\Delta_{+} f_{\alpha}^{-}(t_0) = \Delta_{-} f_{\alpha}^{+}(t_0)$ is monotonically increasing with respect to $\alpha,$ while $\Delta_{-} f_{\alpha}^{-}(t_0) = \Delta_{+} f_{\alpha}^{+}(t_0)$ is monotonically decreasing. In addition, the inequality $\Delta_{+} f_{1}^{-}(t_0) \leq \Delta_{+} f_{1}^{+}(t_0)$. In this setting, we have
$$
\left[\Delta_{gH} f(t_0)\right]_{\alpha} = [\Delta_{+} f_{\alpha}^{-}(t_0), \Delta_{+} f_{\alpha}^{+}(t_0)] = [\Delta_{-} f_{\alpha}^{+}(t_0), \Delta_{-} f_{\alpha}^{-}(t_0)]
$$
for all $\alpha \in [0,1]$.

\item There exists a threshold \( \alpha_0 \in (0,1] \) such that, for all \( \alpha \in [\alpha_0,1] \), the functions  $f_{\alpha }^{-}$ and $f_{\alpha }^{+}$ are uniformly $\Delta$-differentiable at $t_0$   with respect to $\alpha,$ and 
$$
\left[\Delta_{gH} f(t_0)\right]_{\alpha} =\Delta f_{\alpha }^{-}(t_0)=\Delta f_{\alpha }^{+}(t_0)
.$$
\noindent On the other hand, for \( \alpha \in [0,\alpha_0) \), neither \( f_{\alpha}^{-} \) nor \( f_{\alpha}^{+} \) is left or right $\Delta$-differentiable at $t_0$. However, their associated slope functions, \( \varphi_{f_{\alpha}^{-}} \) and \( \varphi_{f_{\alpha}^{+}} \), are  left and right complementary at zero, that is, there exist constants \( a_{\alpha} < b_{\alpha} \) such that
\[
C_{R(0)}(\varphi_{f_{\alpha}^{-}}) = C_{R(0)}(\varphi_{f_{\alpha}^{+}}) = C_{L(0)}(\varphi_{f_{\alpha }^{-}}) = C_{L(0)}(\varphi_{f_{\alpha }^{+}}) = \{a_{\alpha}, b_{\alpha}\},
\]
and the limits 
\[
\lim_{\substack{h \to 0^+\\ h \in \tilde{\mathbb{T}}}} \min\left\{ \varphi_{f_{\alpha}^{-}}(h), \varphi_{f_{\alpha}^{+}}(h) \right\}=\lim_{\substack{h \to 0^-\\ h \in \tilde{\mathbb{T}}}} \min\left\{ \varphi_{f_{\alpha}^{-}}(h), \varphi_{f_{\alpha}^{+}}(h) \right\}=a_{\alpha},\]
\[
\lim_{\substack{h \to 0^+\\ h \in \tilde{\mathbb{T}}}} \max\left\{ \varphi_{f_{\alpha}^{-}}(h), \varphi_{f_{\alpha}^{+}}(h) \right\}=\lim_{\substack{h \to 0^-\\ h \in \tilde{\mathbb{T}}}} \max\left\{ \varphi_{f_{\alpha}^{-}}(h), \varphi_{f_{\alpha}^{+}}(h) \right\}=b_{\alpha} .\]
exist uniformly with respect to  $\alpha$. 

As a result, the  \(\Delta_{gH}\)-derivative of \( f \) at \( t_0 \) is given by
\[
[\Delta_{gH_+} f(t_0)]_{\alpha} = [a_{\alpha}, b_{\alpha}]
\]
for all \( \alpha \in [0,1] \).
\end{enumerate}
\end{theorem}
\begin{theorem}\label{denseder}
Let $f:\mathbb{T} \to \mathbb{R}_{\cal{F}}$ be a fuzzy number-valued function, and let $t_0 \in \mathbb{T}^{\kappa}$ be a  dense point. Then, $f$ is  $\Delta_{gH}$-differentiable at  $t_0$ if one of the following cases holds:
\begin{enumerate}[label=(\roman*)]
\item The functions $f_{\alpha }^{-}$ and $f_{\alpha }^{+}$
are uniformly $\Delta$-differentiable at $t_0$  with respect to $\alpha\in[0,1].$  Furthermore, as functions of  $\alpha$, $\Delta f_{\alpha}^{-}(t_0)$ is monotonically increasing, while $\Delta f_{\alpha}^{+}(t_0)$ is monotonically decreasing. In addition, the
inequality $\Delta f_{1}^{-}(t_0) \leq \Delta f_{1}^{+}(t_0)$ holds. Then, for each \( \alpha \in [0,1] \), the \( \alpha \)-level set of the \( \Delta_{gH} \)-derivative of \( f \) at \( t_0 \) is  given by
\begin{equation*}
    \left[\Delta_{gH} f(t_0)\right]_{\alpha} = \left[\Delta f_{\alpha}^{-}(t_0), \Delta f_{\alpha}^{+}(t_0)\right].
\end{equation*}

\item The functions $f_{\alpha }^{-}$ and $f_{\alpha }^{+}$
are uniformly  $\Delta$-differentiable at $t_0$ with respect to $\alpha\in[0,1].$   Furthermore, as functions of  $\alpha$, $\Delta f_{\alpha}^{+}(t_0)$ is monotonically increasing, while $\Delta f_{\alpha}^{-}(t_0)$ is monotonically decreasing. In addition, the
inequality $\Delta f_{1}^{+}(t_0) \leq \Delta f_{1}^{-}(t_0)$ holds. Then, for each \( \alpha \in [0,1] \), the \( \alpha \)-level set of the \( \Delta_{gH} \)-derivative of \( f \) at \( t_0 \) is  given by
\begin{equation*}
    \left[\Delta_{gH} f(t_0)\right]_{\alpha} = \left[\Delta f_{\alpha}^{+}(t_0), \Delta f_{\alpha}^{-}(t_0)\right].
\end{equation*}

\item  The one sided derivatives  $\Delta_{-} f_{\alpha}^{-}(t_0),$ $\Delta_{+} f_{\alpha}^{-}(t_0)$, $\Delta_{-} f_{\alpha}^{+}(t_0)$, and $\Delta_{+} f_{\alpha}^{+}(t_0)$  exist uniformly with respect to $\alpha\in[0,1]$. Furthermore, the function 
$\Delta_{-} f_{\alpha }^{-}(t_0)=\Delta_{+} f_{\alpha }^{+}(t_0)$ is monotonically increasing with respect to $\alpha,$ while $\Delta_{+} f_{\alpha }^{-}(t_0)=\Delta_{-} f_{\alpha }^{+}(t_0)$ is monotonically decreasing. In addition, the
inequality $\Delta_{-} f_{1}^{-}(t_0)\leq \Delta_{-} f_{1}^{+}(t_0)$ holds. In this setting, we have
$$\left[\Delta_{gH}f(t_0)\right]_{\alpha}=[\Delta_{-} f_{\alpha }^{-}(t_0),\Delta_{-} f_{\alpha }^{+}(t_0)]=[\Delta_{+} f_{\alpha }^{+}(t_0),\Delta_{+} f_{\alpha }^{-}(t_0)]$$
for all $\alpha\in[0,1].$

\item The one sided derivatives  $\Delta_{-} f_{\alpha}^{-}(t_0)$, $\Delta_{+} f_{\alpha}^{-}(t_0)$, $\Delta_{-} f_{\alpha}^{+}(t_0)$, and $\Delta_{+} f_{\alpha}^{+}(t_0)$ exist  uniformly with respect to $\alpha \in [0,1]$.  Furthermore, the function  $\Delta_{+} f_{\alpha}^{-}(t_0) = \Delta_{-} f_{\alpha}^{+}(t_0)$ is monotonically increasing with respect to $\alpha,$ while $\Delta_{-} f_{\alpha}^{-}(t_0) = \Delta_{+} f_{\alpha}^{+}(t_0)$ is monotonically decreasing. In addition, the inequality  $\Delta_{+} f_{1}^{-}(t_0) \leq \Delta_{+} f_{1}^{+}(t_0)$ holds. In this setting, we have
$$
\left[\Delta_{gH} f(t_0)\right]_{\alpha} = [\Delta_{+} f_{\alpha}^{-}(t_0), \Delta_{+} f_{\alpha}^{+}(t_0)] = [\Delta_{-} f_{\alpha}^{+}(t_0), \Delta_{-} f_{\alpha}^{-}(t_0)]
$$
for all $\alpha \in [0,1]$.

\item  There exists a threshold \( \alpha_0 \in (0,1] \) such that, for all \( \alpha \in [\alpha_0,1] \), the functions  $f_{\alpha }^{-}$ and $f_{\alpha }^{+}$ are uniformly $\Delta$-differentiable at $t_0$   with respect to $\alpha,$ and 
$$
\left[\Delta_{gH} f(t_0)\right]_{\alpha} =\Delta f_{\alpha }^{-}(t_0)=\Delta f_{\alpha }^{+}(t_0)
.$$
\noindent On the other hand, for \( \alpha \in [0,\alpha_0) \), neither \( f_{\alpha}^{-} \) nor \( f_{\alpha}^{+} \) is left or right $\Delta$-differentiable at $t_0$. However, their associated slope functions, \( \varphi_{f_{\alpha}^{-}} \) and \( \varphi_{f_{\alpha}^{+}} \), are  left and right complementary at zero, that is, there exist constants \( a_{\alpha} < b_{\alpha} \) such that
\[
C_{R(0)}(\varphi_{f_{\alpha}^{-}}) = C_{R(0)}(\varphi_{f_{\alpha}^{+}}) = C_{L(0)}(\varphi_{f_{\alpha }^{-}}) = C_{L(0)}(\varphi_{f_{\alpha }^{+}}) = \{a_{\alpha}, b_{\alpha}\},
\]
and the limits 
\[
\lim_{\substack{h \to 0^+\\ h \in \tilde{\mathbb{T}}}} \min\left\{ \varphi_{f_{\alpha}^{-}}(h), \varphi_{f_{\alpha}^{+}}(h) \right\}=\lim_{\substack{h \to 0^-\\ h \in \tilde{\mathbb{T}}}} \min\left\{ \varphi_{f_{\alpha}^{-}}(h), \varphi_{f_{\alpha}^{+}}(h) \right\}=a_{\alpha},\]
\[
\lim_{\substack{h \to 0^+\\ h \in \tilde{\mathbb{T}}}} \max\left\{ \varphi_{f_{\alpha}^{-}}(h), \varphi_{f_{\alpha}^{+}}(h) \right\}=\lim_{\substack{h \to 0^-\\ h \in \tilde{\mathbb{T}}}} \max\left\{ \varphi_{f_{\alpha}^{-}}(h), \varphi_{f_{\alpha}^{+}}(h) \right\}=b_{\alpha} .\]
exist uniformly with respect to  $\alpha$. 

\noindent Furthermore, the endpoints \( a_{\alpha} \) and \( b_{\alpha} \) vary monotonically with  \( \alpha \): the function \( a_{\alpha} \) increases, while \( b_{\alpha} \) decreases. These  endpoints also satisfy the consistency condition \( a_1 \leq b_1 \). 

\noindent As a result, the  \(\Delta_{gH}\)-derivative of \( f \) at \( t_0 \) is given by
\[
\left[\Delta_{gH_-}f(t_0)\right]_{\alpha} = [a_{\alpha},b_{\alpha}]
\]
for all $\alpha \in [0,1]$.
\end{enumerate}
\end{theorem}
At any isolated point 
 $t_0\in\mathbb{T}^{\kappa},$  the following theorem is valid.

\begin{theorem}\cite{Mert1}
Let $f:\mathbb{T} \to \mathbb{R}_{\cal{F}}$ be a fuzzy number-valued function, and let $t_0 \in \mathbb{T}^{\kappa}$ be an isolated point. Then, $f$ is  $\Delta_{gH}$-differentiable at  $t_0$, and for each \( \alpha \in [0,1] \), the \( \alpha \)-level set of the \( \Delta_{gH} \)-derivative of \( f \) at \( t_0 \) is an interval given by either
\[
[\Delta_{gH} f(t_0)]_{\alpha} = \left[\Delta f_{\alpha}^{-}(t_0), \Delta f_{\alpha}^{+}(t_0)\right]
\]
or
\[
[\Delta_{gH} f(t_0)]_{\alpha} = \left[\Delta f_{\alpha}^{+}(t_0), \Delta f_{\alpha}^{-}(t_0)\right].
\]
\end{theorem}

%\begin{theorem}\cite{Mert1}
%Let $f:\mathbb{T} \to \mathbb{R}_{\cal{F}}$ be a fuzzy number-valued function, and let $t_0 \in \mathbb{T}^{\kappa}$ be an isolated point. Then, $f$ is  $\Delta_{gH}$-differentiable at  $t_0$ if one of the following two cases holds:
%\begin{enumerate}
%\item The functions $f_{\alpha}^{-}$ and $f_{\alpha}^{+}$ are $\Delta$-differentiable at $t$, uniformly with respect to $\alpha \in [0,1]$. Moreover, $\Delta f_{\alpha}^{-}(t)$ is monotonically increasing and $\Delta f_{\alpha}^{+}(t)$ is monotonically decreasing as functions of $\alpha$, and it holds that $\Delta f_{1}^{-}(t) \leq \Delta f_{1}^{+}(t)$. In this case, we have the following:
%\begin{equation*}
%    \left[\Delta_{gH} f(t)\right]_{\alpha} = \left[\Delta f_{\alpha}^{-}(t), \Delta f_{\alpha}^{+}(t)\right]
%\end{equation*}
%for all $\alpha \in [0,1]$.

%\item The functions $f_{\alpha}^{-}$ and $f_{\alpha}^{+}$ are $\Delta$-differentiable at $t$, uniformly with respect to $\alpha \in [0,1]$. Moreover, $\Delta f_{\alpha}^{+}(t)$ is monotonically increasing and $\Delta f_{\alpha}^{-}(t)$ is monotonically decreasing as functions of $\alpha$, and it holds that  $\Delta f_{1}^{+}(t) \leq \Delta f_{1}^{-}(t)$. In this case, we have  the following:
%\begin{equation*}
%    \left[\Delta_{gH} f(t)\right]_{\alpha} = \left[\Delta f_{\alpha}^{+}(t), \Delta f_{\alpha}^{-}(t)\right]
%\end{equation*}
%for all $\alpha \in [0,1]$.
%\end{enumerate}
%\end{theorem}

\section{Conclusion}

In this paper, we have investigated the generalized Hukuhara  differentiability of fuzzy number-valued functions on arbitrary time scales within the framework of delta calculus. Building upon the foundational contributions of previous works, we have addressed key shortcomings and redundancies. Our main contribution is the development of a unified and comprehensive characterization theorem that captures a wide range of gH-differentiability behaviors, including cases that had previously been overlooked or inadequately described.

The proposed framework resolves inconsistencies related to the treatment of endpoint functions and discontinuities, offering a more flexible and accurate approach to fuzzy differentiability. Importantly, our results are valid across discrete, continuous, and hybrid time scales, making them broadly applicable to both theoretical analysis and practical modeling in dynamic systems involving uncertainty.

This work not only strengthens the mathematical foundations of fuzzy differential calculus on time scales but also opens avenues for further research. Potential directions include the study of higher-order gH-derivatives, the development of numerical methods for fuzzy dynamic equations, and applications in fields such as control theory, economics, and engineering, where fuzzy modeling on diverse time domains is essential.

\section*{\bf{Declarations}}

\noindent {\bf{Funding}}  The authors did not receive any financial support for this work.

\noindent {\bf{Data Availability}} No data has been employed in this article.

\noindent {\bf{Ethical Approval}} Not applicable

\noindent {\bf{Competing Interests}} The author declares no competing interests.


\begin{thebibliography}{99}

\bibitem{Hilger1} S. Hilger, \textit{Ein Makettenkalkuls mit Anwendung auf Zentrumsmannigfaltigkeiten}, Ph.D. Thesis, Universität Würzburg, Würzburg, Germany, 1988.
\bibitem{Peterson} M. Bohner and A. Peterson, \textit{Dynamic Equations on Time Scales: An Introduction with Applications}, Birkhäuser, Boston, 2001.

\bibitem{Peterson1} M. Bohner and A. Peterson, \textit{Advances in Dynamic Equations on Time Scales}, Birkhäuser, Boston, 2002.

\bibitem{Zadeh} L. A. Zadeh, Fuzzy sets, \textit{Information and Control}, \textbf{8} (1965), 338--353.

\bibitem{Klir} G. Klir and B. Yuan, \textit{Fuzzy Sets and Fuzzy Logic}, Vol. 4, Prentice Hall, New Jersey, 1995.

\bibitem{Gomes} L. T. Gomes, L. C. de Barros, and B. Bede, \textit{Fuzzy Differential Equations in Various Approaches}, Springer, Berlin, 2015.

\bibitem{Bedebook} B. Bede, \textit{Mathematics of Fuzzy Sets and Fuzzy Logic}, Springer-Verlag, Berlin, Heidelberg, 2013.

\bibitem{Cano} Y. Chalco-Cano, R. Rodríguez-López, and M. D. Jiménez-Gamero, Characterizations of generalized differentiable fuzzy functions, \textit{Fuzzy Sets and Systems}, \textbf{295} (2016), 37--56.

\bibitem{qiu}
D. Qiu,
Characterizations of generalized Hukuhara differentiability for interval-valued functions,
\textit{Fuzzy Sets and Systems}, \textbf{160}(17), 2009, pp. 2394--2410.

\bibitem{longo} F. Longo, B. Laiate, M. C. Gadotti, J. F. da C.A. Meyer, Characterization results of generalized differentiabilities of fuzzy
functions, Fuzzy Sets and Systems 490 (2024) 109038.


\bibitem{Negoita} C. Negoita and D. Ralescu, \textit{Application of Fuzzy Sets to System Analysis}, Wiley, New York, 1975.

\bibitem{Stef} L. Stefanini, A generalization of Hukuhara difference and division for interval and fuzzy arithmetic, \textit{Fuzzy Sets and Systems}, \textbf{161}(11) (2010), 1564--1584.


\bibitem{Bede} B. Bede and L. Stefanini, Generalized differentiability of fuzzy-valued functions, \textit{Fuzzy Sets and Systems}, \textbf{230} (2013), 119--141.





\bibitem{Diamond2} P. Diamond and P. Kloeden, Metric spaces of fuzzy sets, \textit{Fuzzy Sets and Systems}, \textbf{35}(2) (1990), 241--249.



\bibitem{Zakon} E. Zakon, Mathematical Analysis, Volume I, Trillia Group, West Lafayette, 2004.

\bibitem{Mert1} F. R. Mert, S. Bayeg, B. Kaymakcalan, Generalized  Hukuhara Delta Derivatives: New Characterizations and Extensions for Fuzzy Functions on Time Scales, submitted.







\end{thebibliography}
\end{document}